\documentclass{amsart}

\theoremstyle{plain}
\newtheorem{theorem}{Theorem}[section]
\newtheorem{proposition}{Proposition}[section]
\newtheorem{conjecture}{Conjecture}[section]
\newtheorem{lemma}{Lemma}[section]
\newtheorem{corollary}{Corollary}[section]
\numberwithin{equation}{section}

\newcommand{\qbin}[2]{\genfrac{[}{]}{0pt}{}{#1}{#2}}
\newcommand{\qbins}[2]{{\textstyle\genfrac{[}{]}{0pt}{}{#1}{#2}}}
\newcommand{\Z}{\mathbb{Z}}
\newcommand{\T}{\mathcal{T}}
\newcommand{\U}{\mathcal{U}}
\newcommand{\V}{\mathcal{V}}
\newcommand{\B}{\mathcal{B}}
\newcommand{\Q}{\mathcal{Q}}
\renewcommand{\S}{\mathcal{S}}
\newcommand{\I}{\mathcal{I}}
\newcommand{\cf}{\text{cf}}
\newcommand{\sign}{\text{sign}}

\begin{document}

\title[The generalized Borwein conjecture]{The generalized Borwein
conjecture. II. Refined $q$-trinomial coefficients}

\author[Ole Warnaar]{S. Ole Warnaar}

\address{Department of Mathematics and Statistics,
The University of Melbourne, Vic 3010, Australia}
\email{warnaar@ms.unimelb.edu.au}


\thanks{Work supported by the Australian Research Council}

\begin{abstract}
Transformation formulas for four-parameter refinements of the 
$q$-trinomial coefficients are proven. 
The iterative nature of these transformations allows for the easy
derivation of several infinite series of $q$-trinomial identities,
and can be applied to prove many instances of Bressoud's
generalized Borwein conjecture.
\end{abstract}

\maketitle

\section{Introduction}
This is the second in a series of papers addressing Bressoud's generalized
Borwein conjecture. Defining the Gaussian polynomial or $q$-binomial
coefficient as
\begin{equation}\label{qbin}
\qbin{m+n}{m}=\begin{cases}\displaystyle 
\prod_{k=1}^m\frac{1-q^{n+k}}{1-q^k} & m,n\in\Z_{+} \\
0 & \text{otherwise,}
\end{cases}
\end{equation}
(with $\Z_{+}=\{0,1,2,\dots\}$)
Bressoud \cite{Bressoud96} considered the polynomials
\begin{equation*}
G(N,M;\alpha,\beta,K)=
\sum_{j=-\infty}^{\infty}(-1)^j 
q^{Kj((\alpha+\beta)j+\alpha-\beta)/2}\qbin{M+N}{N-Kj}.
\end{equation*}
Writing $P\geq 0$ if $P$ is a polynomial with nonnegative 
coefficients, he then conjectured the truth of the
following statement concerning $G$.
\begin{conjecture}\label{GBC}
Let $K$ be a positive integer and $N,M,\alpha K,\beta K$
be nonnegative integers such that $1\leq \alpha+\beta\leq 2K-1$ (strict
inequalities when $K=2$) and $\beta-K\leq N-M\leq K-\alpha$.
Then $G(N,M;\alpha,\beta,K)\geq 0$.
\end{conjecture}
This generalizes an earlier conjecture of P.~Borwein \cite{Andrews95}
stating that the coefficients of the polynomials $A_n(q)$, $B_n(q)$ and 
$C_n(q)$, defined by
\begin{equation*}
\prod_{k=1}^n(1-q^{3k-2})(1-q^{3k-1})=A_n(q^3)-qB_n(q^3)-q^2C_n(q^3)
\end{equation*}
are all nonnegative. By the $q$-binomial theorem it readily follows 
that \cite{Andrews95}
\begin{align*}
A_n(q)&=G(n,n;4/3,5/3,3)\\
B_n(q)&=G(n+1,n-1;2/3,7/3,3)\\
C_n(q)&=G(n+1,n-1;1/3,8/3,3).
\end{align*}
For a more comprehensive introduction to the above conjectures we refer
to our first paper in this series \cite{W00b} and to the original
publications by Andrews~\cite{Andrews95} and Bressoud~\cite{Bressoud96}.

Several special cases of the generalized Borwein conjecture have already 
been settled in the literature. When $\alpha$ and $\beta$ are integers the
polynomial $G(N,M;\alpha,\beta,K;q)$ has a combinatorial interpretation
as the generating function of partitions that fit in a rectangle of 
dimensions $M\times N$ and satisfy certain restrictions on their 
hook-differences \cite{ABBBFV87}. For later reference and comparison we 
formalize the $M=N$ case of this in a theorem.
\begin{theorem}\label{T1}
$G(M,M;\alpha,\beta,K)\geq 0$ for $\alpha,\beta,K\in\Z$ such that 
$0\leq \alpha,\beta\leq K$.
\end{theorem}
When at least one of $\alpha$ and $\beta$ is fractional, no combinatorial
interpretation of $G(N,M;\alpha,\beta,K;q)$ is known, except for a few
very simple cases. $G(M,M;1/2,1,2)$, for example, is the generating
function of partitions with largest part at most $M$ and no
parts below its Durfee square.
Despite this lack of a combinatorial interpretation,
Ismail, Kim and Stanton \cite[Thm. 5]{IKS99} 
have proven Conjecture~\ref{GBC} to hold for $\alpha+\beta=K$ with 
$\alpha=(K-N+M\pm 1)/2$ and $M+N$ even.
Again we put the $M=N$ case of this in a theorem.
Because of the symmetry 
\begin{equation}\label{Gsym}
G(M,M;\alpha,\beta,K)=G(M,M;\beta,\alpha,K)
\end{equation}
we may without loss of generality assume $\alpha=(K-1)/2$.
\begin{theorem}\label{T2}
For $K$ a positive integer
$G(M,M;(K-1)/2,(K+1)/2,K)\geq 0$.
\end{theorem}
When $K$ is odd this is of course contained in Theorem~\ref{T1}.
Finally we quote a result obtained in our first paper by use of the 
Burge transform \cite[Cor. 3.2]{W00b}.
\begin{theorem}\label{T3}
$G(M,M;b,b+1/a,a)\geq 0$ for $a,b$ coprime integers such that $0<b<a$.
\end{theorem}
Similar results were obtained for $N\neq M$ with
both $\alpha$ and $\beta$ noninteger \cite[Cor 5.1]{W00b}. 
It is quite clear, however, that proving $G(M,M;\alpha,\beta,K)\geq 0$
when both $\alpha$ and $\beta$ are fractional and not $\alpha=(K-1)/2$
and $\beta=(K+1)/2$ --- $A_n(q)$ of the original Borwein conjecture
falls in this class --- is rather more difficult.
In this paper new transformation formulas will be applied to 
make some progress in this direction.
To state our results we define
\begin{subequations}\label{dmu}
\begin{align}
d_{\bar{a}}=d&=a_0+a_1+\cdots+a_n  \\
\mu_{\bar{a}}=\mu&=\frac{2}{3}(4^{d-1}-1)-\sum_{i=1}^n 4^{a_i+\cdots+a_n-2},
\end{align}
\end{subequations}
where $\bar{a}=(a_0,a_1,\dots,a_n)\in\Z_{+}^{n+1}$.
Note that $\mu\in\Z_{+}$ provided $a_n\geq 2-\delta_{n,0}$.
With these definitions our main results are the following three theorems,
which generalize Theorems~\ref{T1}--\ref{T3}.
\begin{theorem}\label{T4}
For $n\geq 0$, let $\bar{a}=(a_0,\dots,a_n)\in\Z_{+}^{n+1}$ such that
$a_0\geq 0$, $a_1,\dots,a_{n-1}\geq 1$ and $a_n\geq 2-\delta_{n,0}$. 
Then 
\begin{equation*}
G(M,M;(\alpha+\mu K)/2^{d-1},(\beta+\mu K)/2^{d-1},2^{d-1}K)\geq 0,
\end{equation*}
for $\alpha,\beta,K\in\Z$ such that $0\leq \alpha,\beta\leq K$,
and $d$ and $\mu$ given by \eqref{dmu}.
\end{theorem}
For $\bar{a}=(1)$ there holds $d=1$ and $\mu=0$ so that we recover
Theorem~\ref{T1}.
\begin{theorem}\label{T5}
With the same conditions as in Theorem~\ref{T4} there holds
\begin{equation*}
G(M,M;((2\mu+1)K-1)/2^d,((2\mu+1)K+1)/2^d,2^{d-1}K)\geq 0,
\end{equation*}
for $K$ a positive integer and $d$ and $\mu$ given by \eqref{dmu}.
\end{theorem}
For $\bar{a}=(1)$ this reduces to Theorem~\ref{T2}.
\begin{theorem}\label{T6}
With the same conditions as in Theorem~\ref{T4} there holds
\begin{equation*}
G(M,M;b/2^{d-1},(b+1/a)/2^{d-1},2^{d-1}a)\geq 0,
\end{equation*}
for $a,b$ coprime integers such that $\mu a<b<(\mu+1)a$,
and $d$ and $\mu$ given by \eqref{dmu}.
\end{theorem}
For $\bar{a}=(1)$ this reduces to Theorem~\ref{T3}.
A slight reformulation of Theorem~\ref{T6} will be given
in Theorem~\ref{T8} of section~\ref{sec6}.

\subsection*{Outline of the paper}
The first part of this paper deals with the theory of $q$-trinomial 
coefficients. In the next section we review the basics of 
$q$-trinomial coefficients and then extend the 
theory to refined $q$-trinomial coefficients. 
Our main results are Theorems~\ref{mainthm} and \ref{mainthm2}, 
which are two elegant transformation formulas for refined $q$-trinomial 
coefficients that can be viewed as trinomial analogues of the Burge 
transform.
The technical section \ref{sec3} contains proofs of 
some of our main claims concerning refined $q$-trinomials.

The second part of the paper contains applications of
the transformation formulas of section~\ref{sec2}, with
sections~\ref{sec4} and \ref{sec5} devoted to proving $q$-trinomial 
identities and Rogers--Ramanujan-type identities, and sections~\ref{sec6}
and \ref{sec7} devoted to the generalized Borwein conjecture.
In the appendix some simple summation formulae needed in the main text 
are established.

\section{Refined $q$-Trinomial coefficients}\label{sec2}
We employ the following standard notations for the $q$-shifted factorial:
$(a;q)_n=(a)_n=\prod_{j=1}^n (1-aq^{j-1})$ for $n\geq 0$,
$(a;q)_n=(a)_n=1/(aq^n;q)_{-n}$ for $n\in\Z$ (so that $1/(q)_{-n}=0$
for $n>0$) and $(a_1,\dots,a_k;q)_n=(a_1,\dots,a_k)_n=(a_1)_n\cdots (a_k)_n$.
Whenever series are nonterminating series
it is tacitly assumed that $|q|<1$.

In analogy with the definition of binomial coefficients,
the trinomial coefficients $\binom{L}{a}_2$ are defined by the expansion
\begin{equation}\label{trinexp}
(1+x+x^2)^L=\sum_{a=-L}^L\binom{L}{a}_2 x^{a+L}.
\end{equation}
Double application of the binomial expansion shows that
\begin{equation}\label{trinomials}
\binom{L}{a}_2=\sum_{k=0}^L \binom{L}{k}\binom{L-k}{k+a}.
\end{equation}
The analogy with binomials breaks down when it comes to defining
$q$-analogues. The binomial expansion is readily 
generalized to the $q$-case by \cite[(3.3.6)]{Andrews76}
\begin{equation}\label{qbt}
(x)_L=\sum_{a=0}^L (-x)^a q^{\binom{a}{2}}\qbin{L}{a},
\end{equation}
but no $q$-analogue of \eqref{trinomials} seems possible that yields a 
$q$-version of \eqref{trinexp}. Despite this complication, Andrews
and Baxter \cite{AB87} successfully defined useful $q$-trinomial 
coefficients. Here we need just two of the simplest $q$-analogues
of \eqref{trinomials} given by \cite[Eq. (2.7); $B=A$]{AB87}
\begin{equation}\label{qt2}
\qbin{L}{a}_{2;q}=\qbin{L}{a}_2
=\sum_{k=0}^L q^{k(k+a)}\qbin{L}{k}\qbin{L-k}{k+a}
\end{equation}
and \cite[Eq. (2.8)]{AB87}
\begin{equation}\label{Td}
T(L,a;q)=T(L,a)=q^{\frac{1}{2}(L^2-a^2)}\qbin{L}{a}_{2;q^{-1}}.
\end{equation}
An explicit expression for $T(L,a)$ needed later is given by 
\cite[Eq. (2.60)]{AB87}
\begin{equation}\label{qtT}
T(L,a)=\sum_{\substack{n=0 \\n+a+L \text{ even}}}^L
q^{\frac{1}{2}n^2}\qbin{L-n}{\frac{1}{2}(L-a-n)}\qbin{L}{n}.
\end{equation}

It is easy to see from \eqref{qt2} and \eqref{qtT} that the
$q$-trinomial coefficients obey the symmetry
$\qbins{L}{a}_2=\qbins{L}{-a}_2$ and $T(L,a)=T(L,-a)$.
Almost as easy to establish are the large $L$ limits.
By a limit of the $q$-Gauss sum \cite[Eq. (II.8)]{GR90},
\begin{equation}\label{tlim}
\lim_{L\to\infty}\qbin{L}{a}_2=
\sum_{k=0}^{\infty}\frac{q^{k(k+a)}}{(q)_k(q)_{k+a}}
=\frac{1}{(q)_{\infty}},
\end{equation}
and by Euler's $q$-exponential sum \cite[Eq. (II.2)]{GR90},
\begin{equation*}
\lim_{\substack{L\to\infty \\ L+a+\sigma \text{ even}}}T(L,a)=
\sum_{\substack{n=0 \\ n+\sigma \text{ even}}}^{\infty}
\frac{q^{\frac{1}{2}n^2}}{(q)_n}=
\frac{(-q^{1/2})_{\infty}+(-1)^{\sigma}(q^{1/2})_{\infty}}
{2(q)_{\infty}}.
\end{equation*}

To conclude our brief review of $q$-trinomial coefficients we mention
that in identities one often encounters the same linear combination of 
two such coefficients. For this reason it is helpful to 
define \cite{Andrews90b}
\begin{equation}\label{U}
U(L,a)=T(L,a)+T(L,a+1),
\end{equation}
which has a limiting behaviour somewhat simpler to that of $T(L,a)$,
\begin{equation}\label{Ulim}
\lim_{L\to\infty}U(L,a)=\frac{(-q^{1/2})_{\infty}}{(q)_{\infty}}.
\end{equation}

In the following we go well-beyond $q$-trinomial coefficients, and introduce 
two polynomials $\S$ and $\T$ that can be viewed as four-parameter 
extensions of $\qbins{L}{a}_2$ and $T(L,a)$, respectively.
Assuming that $L,M,a,b$ are all integers we define
\begin{equation}\label{Sdef}
\S(L,M,a,b;q)=\S(L,M,a,b)=\sum_{k=0}^L
q^{k(k+a)}\qbins{L+M-a-2k}{M}\qbins{M-a+b}{k}\qbins{M+a-b}{k+a}
\end{equation}
and
\begin{align}\label{Tdef}
\T(L,M,a,b;q)&=\T(L,M,a,b) \\[2mm]
&=\sum_{\substack{n=0 \\n+a+L \text{ even}}}^L
q^{\frac{1}{2}n^2}
\qbins{M}{n}\qbins{M+b+(L-a-n)/2}{M+b}\qbins{M-b+(L+a-n)/2}{M-b}.
\notag
\end{align}
Comparison with \eqref{qt2} and \eqref{qtT} shows that
\begin{align}\label{SlimM}
(q)_L \lim_{M\to\infty}\S(L,M,a,b)&=\qbin{L}{a}_2
\intertext{and}
(q)_L \lim_{M\to\infty}\T(L,M,a,b)&=T(L,a).\label{TlimM}
\end{align}

Before we list the most important properties of $\S$ and
$\T$ let us remark that the polynomial $\T$ was recently introduced 
in \cite{W00a}. Following the terminology of \cite{W00a} we will call 
$\T(L,M,a,b)$ a refined $q$-trinomial coefficient (for reasons that
will become clear shortly, and \textit{not} because of \eqref{TlimM}).

The first three properties of $\S$ and $\T$ listed below follow 
directly from the definitions. With $\Q$ to mean either $\S$ or $\T$ we 
have the range of support
\begin{subequations}
\begin{align}\label{rangeS}
\S(L,M,a,b)&\neq 0 \;\text{ iff $\;|a|\leq L$, $|b|\leq M$ and 
$|a-b|\leq M$} \\
\label{rangeT}
\T(L,M,a,b)&\neq 0 \;\text{ iff $\;|a|\leq L$, $|b|\leq M$ and 
$\frac{1}{2}(a+L)\in\Z$ if $M=0$,}
\end{align}
\end{subequations}
the symmetry
\begin{equation}\label{Qsymm}
\Q(L,M,a,b)=\Q(L,M,-a,-b)
\end{equation}
and the duality
\begin{equation}\label{Qdual}
\Q(L,M,a,b;1/q)=q^{ab-LM}\Q(L,M,a,b;q).
\end{equation}
Whereas the first two equations are the obvious 
generalizations of analogous properties of $q$-trinomial coefficients,
the last equation is in clear contrast with \eqref{Td}.
The next result, to be compared with \eqref{tlim}, will be important when
we address the generalized Borwein conjecture;
\begin{equation}\label{SlimL}
(q)_M \lim_{L\to\infty}\S(L,M,a,b)=
\sum_{k=0}^{\infty}q^{k(k+a)}\qbins{M-a+b}{k}\qbins{M+a-b}{k+a}
=\qbins{2M}{M-b}.
\end{equation}
Here the second equality follows from the $q$-Chu--Vandermonde sum 
\cite[(II.7)]{GR90}
\begin{equation}\label{qCV}
\sum_{k=0}^n \frac{(a,q^{-n})_k}{(q,c)_k}\Bigl(\frac{cq^n}{a}\Bigr)^k
=\frac{(c/a)_n}{(c)_n},
\end{equation}
with $a\to q^{-(M-a+b)}$, $n\to M-b$ and $c\to q^{a+1}$.

We now come to the main results of this section.
\begin{theorem}\label{mainthm}
For $L,M,a,b\in\Z$ such that $ab\geq 0$
\begin{equation}\label{Tinv}
\sum_{i=0}^M q^{\frac{1}{2}i^2}\qbin{L+M-i}{L}\T(L-i,i,a,b)
=q^{\frac{1}{2}b^2}\T(L,M,a+b,b).
\end{equation}
\end{theorem}
This transformation was announced in \cite[Thm. 3.1]{W00a}.
Our next transform has not appeared before.
\begin{theorem}\label{mainthm2}
For $L,M,a,b\in\Z$ such that $ab\geq 0$, and such that $|a|\leq M$ 
if $|b|\leq M$ and $|a+b|\leq L$, 
\begin{equation}\label{TtoS}
\sum_{i=0}^M q^{\frac{1}{2}i^2}
\qbin{L+M-i}{L}\T(i,L-i,b,a)
= q^{\frac{1}{2}b^2}\S(L,M,a+b,b).
\end{equation}
\end{theorem}
A discussion of the conditions imposed on the parameters 
(which are \textit{not} sharp) precedes the proofs given in sections
\ref{sec31} and \ref{sec32}.

Theorems~\ref{mainthm} and \ref{mainthm2} justify calling $\T(L,M,a,b)$
a refined $q$-trinomial coefficient because they imply
\begin{equation*}
\sum_{i=0}^L q^{\frac{1}{2}(i^2-b^2)}\T(L-i,i,a-b,b)=T(L,a)
\end{equation*}
and
\begin{equation*}
\sum_{i=0}^L q^{\frac{1}{2}(i^2-b^2)} \T(i,L-i,b,a-b)=\qbin{L}{a}_2.
\end{equation*}
Here it is assumed in both formulae that $0\leq b\leq a$ or $a\leq b\leq 0$.
The first equation follows by taking the large $M$ limit in 
\eqref{Tinv} using \eqref{TlimM}. The second equation follows from the first
by application of \eqref{Td} and \eqref{Qdual} or from \eqref{TtoS} by
taking $M$ to infinity and using \eqref{SlimM}.

Given the above two theorems, an obvious question is whether there also
exist transformations from $\S$ to $\T$ or from $\S$ to $\S$. The only 
result we found in this direction is the following not-so-useful summation.
\begin{lemma}
For $L,M,a,b\in\Z$ such that 
$L\leq M$ and $|a-b|\leq L$ if $|a|\leq L$,
$|b|\leq M$ and $|a-b|\leq \max\{L,M\}$,
\begin{equation*}
\sum_{i=0}^{M-L}q^{Li}
\qbin{M-L}{i}\S(L-i,M-i,a,b)=\S(M,L,b,a).
\end{equation*}
\end{lemma}
Since we will not use this transformation we omit its proof.
We note however that for $a=b=0$ it coincides with Theorem~\ref{mainthm3}
below (with $a=b=0$ and $L$ and $M$ interchanged), 
and the proof for more general $a$ and $b$ is a simple modification of 
the proof of that theorem as given in section~\ref{sec33}.

Before we can state our next two results we first need to define
\begin{equation}\label{Bdef}
\B(L,M,a,b;q)=\B(L,M,a,b)=\qbin{M+b+L-a}{M+b}\qbin{M-b+L+a}{M-b},
\end{equation}
for $L,M,a,b\in\Z/2$ such that $L+a$ and $M+b$ are integers.
Note that $\B(L,M,a,b)$ is nonzero for $|b|\leq M$ and $|a|\leq L$ only.

\begin{lemma}
For $L,M,a,b\in\Z$
\begin{equation}\label{BtoT}
\sum_{\substack{n=0 \\n+L+a \text{ \rm even}}}^M
q^{\frac{1}{2}n^2}\qbin{M}{n}\B((L-n)/2,M,a/2,b)=\T(L,M,a,b).
\end{equation}
\end{lemma} 
\begin{proof}
Substituting the definitions of $\T$ and $\B$ gives the desired result.
\end{proof}
\begin{theorem}\label{mainthm3}
For $L,M,a,b\in\Z$ such that 
$M\leq L$ and $|a-b|\leq M$ if $|a|\leq M/2$,
$|a+b|\leq L$ and $|a-b|\leq \max\{L,M\}$, 
\begin{multline}\label{StoS}
\sum_{i=0}^{L-M}\sum_{k=0}^L q^{M i+k^2}
\qbin{L-M}{i}\qbin{L+M-2i-2k}{L-i}\B(k,L-i-k,a,b) \\
=q^{a^2}\S(L,M,a+b,2a).
\end{multline}
\end{theorem}
The conditions imposed on the above summation formula are sharp.
Their origin will be discussed in the proof given in section~\ref{sec33}.

Next we derive two corollaries of Theorems~\ref{mainthm2} 
and \ref{mainthm3}.
First, taking \eqref{TtoS}, inserting the definition of $\T$ and letting
$L$ tend to infinity using \eqref{SlimL} yields
\begin{equation*}
\sum_{i=0}^M 
\sum_{\substack{n=0 \\n+b+i \text{ even}}}^i
q^{\frac{1}{2}(i^2+n^2)}\qbin{M}{i}
\qbin{i-n}{\frac{1}{2}(i-b-n)}\qbin{i}{n}
=q^{\frac{1}{2}b^2}\qbin{2M}{M-b}.
\end{equation*}
Replacing $b$ by $2a$ and $n$ by $i-2k$ gives rise to the following result.
\begin{corollary}\label{cor1}
For $M,a\in\Z$
\begin{equation*}
\sum_{k=0}^{\infty}C_{M,k}(q)\qbin{2k}{k-a}=q^{2a^2}\qbin{2M}{M-2a},
\end{equation*}
where
\begin{equation*}
C_{M,k}(q)=\sum_{i=0}^M q^{(i-k)^2+k^2}\qbin{M}{i}\qbin{i}{2k}\geq 0.
\end{equation*}
\end{corollary}
Taking \eqref{StoS}, inserting the definition of $\B$ and sending
$L$ to infinity using \eqref{SlimL} yields a very similar result.
\begin{corollary}\label{cor2}
For $M,a\in\Z$
\begin{equation*}
\sum_{k=0}^{\infty}
\bar{C}_{M,k}(q)\qbin{2k}{k-a}=q^{a^2}\qbin{2M}{M-2a},
\end{equation*}
where
\begin{equation*}
\bar{C}_{M,k}(q)=\sum_{i=0}^M q^{M(M-i)+k^2}\qbin{M}{i}\qbin{i}{2k}\geq 0.
\end{equation*}
\end{corollary}
This can also be obtained from Corollary~\ref{cor1} by the substitution
$q\to 1/q$.

We conclude our discussion of refined $q$-trinomial coefficients
by introducing the refined version of the polynomial $U$ of 
equation \eqref{U} and another polynomial frequently needed;
\begin{align}\label{Udef}
\U(L,M,a,b)&=\T(L,M,a,b)+\T(L,M,a+1,b) \\
\V(L,M,a,b)&=\S(L,M,a,b)+q^{b+1/2}\S(L,M,a+1,b+1).
\label{Vdef}
\end{align}
The following limits of $\U$ and $\V$ will be useful later
\begin{equation}\label{UlimL}
\lim_{L\to\infty}\U(L,M,a,b)=\frac{(-q^{1/2})_M}{(q)_{2M}}
\qbin{2M}{M-b}
\end{equation}
and
\begin{align}\label{UlimLM}
\lim_{L,M \to \infty}\U(L,M,a,b)
&=\frac{(-q^{1/2})_{\infty}}{(q)^2_{\infty}} \\
\label{VlimLM}
\lim_{L,M \to \infty}\V(L,M,a,b)
&=\frac{1+q^{b+1/2}}{(q)^2_{\infty}}.
\end{align}
The first limit follows from the definitions of $\U$ and $\T$ and the
$q$-binomial theorem \eqref{qbt} with $x=-q^{1/2}$.
The second limit is obvious from the first, 
and the last limit follows from \eqref{tlim} and \eqref{SlimM}.

We finally compare some of our results for refined $q$-trinomial 
coefficients with known results for the polynomial $\B$ of equation 
\eqref{Bdef}. First we note that $\B$ obeys
\begin{align}
\B(L,M,-a,-b)&=\B(L,M,a,b)  \notag \\
\B(L,M,a,b)&=\B(M,L,b,a)  \label{Bsymm} \\
\B(L,M,a,b;1/q)& =q^{2ab-2LM}\B(L,M,a,b;q). \notag
\end{align}
The first and last of these relations are similar to \eqref{Qsymm} and
\eqref{Qdual} satisfied by $\S$ and $\T$. Surprisingly, the analogy goes 
much further, and the following two 
theorems are clear analogues of Theorems~\ref{mainthm} and \ref{mainthm2}.
\begin{theorem}\label{thmBurge}
For $L,M,a,b\in\Z$ such that $|a-b|\leq L$ 
if $|b|\leq M$ and $|a+b|\leq L$, 
\begin{equation*}
\sum_{i=0}^M q^{i^2}\qbin{2L+M-i}{2L}\B(L-i,i,a,b)=q^{b^2}\B(L,M,a+b,b).
\end{equation*}
\end{theorem}
\begin{theorem}\label{thmBurge2}
With the same conditions as above
\begin{equation*}
\sum_{i=0}^M q^{i^2}\qbin{2L+M-i}{2L}\B(i,L-i,b,a)=q^{b^2}\B(L,M,a+b,b).
\end{equation*}
\end{theorem}
These two theorems are known as the Burge transform, see e.g.,
\cite{Burge93,FLW97,SW00,W00b,AB01}. A proof follows from the 
$q$-Saalsch\"utz sum \eqref{qSaal}. Unlike \eqref{Tinv} and 
\eqref{TtoS}, which are independent transformations,
the above transformations imply one another thanks to 
the symmetry \eqref{Bsymm}.
Another more important difference between the Burge transform and the
transformations for $\T$ and $\S$ is that the two Burge transformations 
can be iterated to yield a binary tree of transformations 
\cite{Burge93,FLW97,W00b}, while the transformations 
for $\T$ and $\S$ only give rise to an infinite double chain.

Later in the paper we also need the following extension of 
Theorem~\ref{thmBurge} involving the polynomial
\begin{equation}\label{Brs}
\B_{r,s}(L,M,a,b)=\qbin{M+b+s+L-a+r}{M+b+s}\qbin{M-b+L+a}{M-b}.
\end{equation}
\begin{theorem}[\cite{Burge93,SW00}]\label{thmBurge3}
For $L,M,a,b,r,s\in\Z$ such that $-L\leq a-b-s\leq L+r$
if $-M-s\leq b\leq M$ and $-L\leq a+b \leq L+r$,
\begin{equation*}
\sum_{i=b}^Mq^{i(i+s)}\qbin{2L+M+r-i}{2L+r}\B_{r+s,s}(L-i-s,i,a,b)=
q^{b(b+s)}\B_{r,s}(L,M,a+b,b).
\end{equation*}
\end{theorem}

\section{Proofs of Theorems \ref{mainthm}--\ref{mainthm3}}\label{sec3}
\subsection{Proof of Theorem \ref{mainthm}}\label{sec31}
Before we commence with the proof a few comments are in order.
From equation \eqref{rangeT} it follows that the summand on the left 
is zero if not $|b|\leq i \leq \min\{L-|a|,M\}$.
This implies that the left-hand side vanishes trivially if $|b|>M$ or
$|a|+|b|>L$. By the same equation \eqref{rangeT} the right-hand side
vanishes if $|b|>M$ or $|a+b|>L$.
One might thus hope that a sufficient condition for Theorem~\ref{mainthm}
to hold nontrivially 
would be $|a|+|b|\leq L$ and $|b|\leq M$.
However this appears not to be the case and 
$|a|+|b|\leq L$ needs to be replaced by
(i) $|a|+|b|\leq L$ with $ab\geq 0$,
or (ii) $|a|+2|b|\leq L$ with $ab\leq 0$.
Since in all interesting applications of the theorem it turns out that
$a$ and $b$ have the same signature, we have omitted the
cases where $a$ and $b$ have opposite sign in the statement of the theorem
and in the proof given below.
However, in \cite[Thm. 1.2]{SW00} a very general transformation formula 
is proven which for $N=2$, $M\to M-b$, $\ell\to 2b$, $L_1\to (L-a+2b)/2$
and $L_2\to (L+a-2b)$ coincides with Theorem~\ref{mainthm}.
The condition $L_1,L_2\geq 0$ as given in \cite{SW00}
establishes the validity of \eqref{Tinv} for $|a-2b|\leq L$.

\begin{proof}[Proof of Theorem \ref{mainthm}]
Substituting the definition of $\T$ in \eqref{Tinv} we are to prove
\begin{multline}\label{inn}
\sum_{i,n}q^{\frac{1}{2}(i^2+n^2)}\qbins{L+M-i}{L}\qbins{i}{n}
\qbins{i+b+(L-i-a-n)/2}{i+b}\qbins{i-b+(L-i+a-n)/2}{i-b} \\
=q^{\frac{1}{2}b^2}\sum_{n}q^{\frac{1}{2}n^2}\qbins{M}{n}
\qbins{M+b+(L-a-b-n)/2}{M+b}\qbins{M-b+(L+a+b-n)/2}{M-b},
\end{multline}
with $a,b$ in the ranges specified by the theorem, and where
we assume that $L+i+n+b$ is even on the left and
$L+n+a+b$ is even on the right.

Since both sides of \eqref{inn} are symmetric under
simultaneous negation of $a$ and $b$ we may without loss of generality
assume that $a,b\geq 0$ in the following. In view of the
previous discussion we may also assume that $b\leq M$ and $L\geq 0$.
As a first step we use the symmetry 
\begin{equation}\label{qbinsym}
\qbins{m+n}{m}=\qbins{m+n}{n}
\end{equation}
to rewrite \eqref{inn} in the less-symmetric form 
\begin{multline}\label{inn2}
\sum_{i,n}q^{\frac{1}{2}(i^2+n^2)}\qbins{L+M-i}{L}\qbins{i}{n}
\qbins{i+b+(L-i-a-n)/2}{(L-i-a-n)/2}\qbins{i-b+(L-i+a-n)/2}{i-b} \\
= q^{\frac{1}{2}b^2}\sum_{n}q^{\frac{1}{2}n^2}\qbins{M}{n}
\qbins{M+b+(L-a-b-n)/2}{(L-a-b-n)/2}\qbins{M-b+(L+a+b-n)/2}{(L+a+b-n)/2}.
\end{multline}
According to the definition \eqref{qbin} of the $q$-binomial coefficient,
$\qbins{m+n}{m}$ is zero if $m<0$.
We will now show that in the case of equation \eqref{inn2}
this condition together with $a,b\geq 0$, $b\leq M$ and $L\geq 0$
implies that all of the top-entries
of the various $q$-binomials are nonnegative.
First consider the left side. The summand vanishes if
not both $L-i-a-n$ and $i-b$ are nonnegative.
This gives the following inequalities for the
top-entries of the four $q$-binomials:
$L+M-i\geq L+M-b\geq 0$, $i\geq b\geq 0$, $i+b+(L-i-a-n)/2\geq 2b\geq 0$ and
$i-b+(L-i+a-n)/2\geq a\geq 0$.
For the right side of \eqref{inn2} it is equally simple.
Since $L-a-b-n\geq 0$ one has
$M\geq 0$, $M+b+(L-a-b-n)/2\geq M+b\geq 0$ and
$M-b+(L+a+b-n)/2\geq M+a\geq 0$.

Now recall the following modified definition of the $q$-binomial
coefficient:
\begin{equation}\label{qbinmod}
\qbin{m+n}{m}=
\begin{cases}\displaystyle \frac{(q^{n+1})_m}{(q)_m} &
\text{for $m\in\Z_+$, $n\in\Z$} \\[3mm]
0 & \text{otherwise.} \end{cases}
\end{equation}
The only difference between \eqref{qbin} and \eqref{qbinmod} is that in 
the latter $\qbins{m+n}{m}$ is nonzero for $m\geq 0$ and $m+n<0$.
Since we have just argued that in equation \eqref{inn2} the top-entries
of the $q$-binomials cannot be negative by the conditions on the lower 
entries we may in the remainder of our proof of \eqref{inn2} assume 
definition \eqref{qbinmod}.

After these preliminaries we shall transform the left side of
\eqref{inn2} into the right side. First we make the simultaneous
changes $i\to i+n$ and $n\to i$ to get
\begin{equation}\label{LHS1}
\text{LHS}\eqref{inn2}=
\sum_{i,n}q^{\frac{1}{2}n^2+i(i+n)}\qbins{L+M-n-i}{L}\qbins{i+n}{i}
\qbins{(L-a+n)/2+b}{(L-a-n)/2-i}\qbins{(L+a+n)/2-b}{i+n-b}.
\end{equation}
To proceed we need Sears' ${_4\phi_3}$ transformation 
\cite[Eq. (III.15)]{GR90}
\begin{equation}\label{Searstrafo}
\sum_{k=0}^n\frac{(a,b,c,q^{-n})_k \, q^k}{(q,d,e,f)_k}=
a^n \, \frac{(e/a,f/a)_n}{(e,f)_n} \,
\sum_{k=0}^n \frac{(a,d/b,d/c,q^{-n})_k \, q^k}
{(q,d,aq^{1-n}/e,aq^{1-n}/f)_k},
\end{equation}
true for $def=abcdq^{1-n}$.
Making the substitutions $n\to d-e$, $a\to q^{-c}$, $b\to q^{a+1}$,
$c\to q^{g-f}$, $d\to q^{g+1}$, $e\to q^{-b}$ and $f\to q^{e+1}$ 
results in the following transformation for sums of products of $q$-binomial
coefficients:
\begin{multline}\label{sears}
\sum_{i=0}^{d-e}q^{i(i-c+e+g)}\qbins{a-i}{a-b}\qbins{i+c}{i}
\qbins{d}{d-e-i}\qbins{f}{i+g}\\
=\sum_{i=0}^{d-e}q^{i(i-c+e+g)}\qbins{a-d+e}{i-c+e+f}\qbins{c-g}{i}
\qbins{b+d-i}{d-e-i}\qbins{i+f}{i+g},
\end{multline}
for $a,b,c,d,e,f,g\in\Z$ such that $a+c=b+d+f$.
Returning to \eqref{LHS1}, we utilize \eqref{sears}
to transform the sum over $i$. As a result
\begin{equation*}
\text{LHS}\eqref{inn2}=
\sum_{i,n}q^{\frac{1}{2}n^2+i(i+n)}
\qbins{M+(L+a-n)/2}{i+(L+a+n)/2}\qbins{b}{i}
\qbins{M+b+(L-a-n)/2-i}{(L-a-n)/2-i}\qbins{i-b+(L+a+n)/2}{i+n-b}.
\end{equation*}
By the simultaneous changes $i\to n-i$ and $n\to 2i-n+b$ this yields
\begin{multline*}
\text{LHS}\eqref{inn2}=q^{\frac{1}{2}b^2}
\sum_n q^{\frac{1}{2}n^2}
\qbins{M+b+(L-a-b-n)/2}{(L-a-b-n)/2} \\
\times \sum_i q^{i(i-n+b)}
\qbins{M+(L+a-b+n)/2-i}{(L+a+b+n)/2}\qbins{b}{n-i}\qbins{(L+a-b+n)/2}{i}.
\end{multline*}
Now set $c=g=0$ in \eqref{sears} and eliminate $b$.
This yields the well-known \cite{Gould72,Carlitz74,Andrews76}
\begin{equation}\label{qS}
\sum_{i=0}^{d-e}q^{i(i+e)}\qbins{a-i}{d+f}\qbins{d}{d-e-i}\qbins{f}{i}
=\qbins{a-d+e}{e+f}\qbins{a-f}{d-e},
\end{equation}
which, because of definition \eqref{qbinmod}, is valid for all $a,d,e,f\in\Z$.
Using \eqref{qS} to carry out the sum over $i$ results in the right-hand
side of \eqref{inn2}.
\end{proof}

\subsection{Proof of Theorem \ref{mainthm2}}\label{sec32}
The first part of the discussion at the start of section~\ref{sec31}
also applies here.
That is, the content of Theorem~\ref{mainthm2} is nontrivial for
$|b|\leq M$ and $|a|+|b|\leq L$ only.
It is also true again that the bounds are not sharp, since for
some $a$ and $b$ of opposite sign, such that
$|a|+|b|\leq L$, the theorem holds as well.
The extra condition (compared to Theorem~\ref{mainthm}) $|a|\leq M$
when both $|b|\leq M$ and $|a|+|b|\leq L$ does appear to be sharp.

\begin{proof}
Without loss of generality we may assume $a,b\geq 0$, $b\leq M$ and
$L\geq 0$. Of course we also have $a\leq M$ as a necessary condition.
Substituting the definitions of $\T$ and $\S$ in the above identity
and using \eqref{qbinsym} to asymmetrize, we are to prove
\begin{multline}\label{ink}
\sum_{i,n} q^{\frac{1}{2}(i^2+n^2)}
\qbins{L+M-i}{L} \qbins{L-i}{n}
\qbins{L-i+a+(i-b-n)/2}{(i-b-n)/2}\qbins{L-i-a+(i+b-n)/2}{L-i-a} \\
=q^{\frac{1}{2}b^2}\sum_k
q^{k(k+a+b)}\qbins{L+M-a-b-2k}{L-a-b-2k}\qbins{M-a}{k}\qbins{M+a}{k+a+b},
\end{multline}
where the parity rule that $i+n+b$ on the left must be even is implicit.

As in the proof of Theorem~\ref{mainthm} we will now show that
the top-entries of all seven $q$-binomial coefficients are
nonnegative by the conditions on the parameters and by the
condition that the lower entries are nonnegative.
This allows us to again assume the modified definition \eqref{qbinmod}
of the $q$-binomials in our proof.

First consider the left-hand side. Because $(i-b-n)/2\geq 0$ and
$L-i-a\geq 0$ we get $L+M-i\geq M+a\geq 0$,
$L-i\geq a\geq 0$, $L-i+a+(i-b-n)/2\geq 2a\geq 0$ and
$L-i-a+(i+b-n)/2 \geq b$.
On the right-hand side, since $L-a-b-2k\geq 0$, we get
$L+M-a-b-2k\geq M$ for the top-entry of the first $q$-binomial.
It is the above-discussed extra condition $a\leq M$ that ensures that also
the second $q$-binomial on the right has a nonnegative top-entry.
Since clearly also $M+a$ is nonnegative we are indeed in a position to
assume \eqref{qbinmod} in the remainder of the proof.

We begin with the simultaneous changes $i\to i+k+b$ and $n\to k-i$ to find
\begin{equation*}
\text{LHS}\eqref{ink}
=\sum_{i,k}q^{\frac{1}{2}b^2+k(k+b)+i(i+b)}
\qbins{L+M-b-k-i}{L}\qbins{L-b-k-i}{k-i}
\qbins{L+a-b-k}{i}\qbins{L-a-k}{L-a-b-k-i}.
\end{equation*}
Next we apply the following transformation formula similar to \eqref{sears},
which can again be viewed as a corollary of the Sears transform
\eqref{Searstrafo}:
\begin{equation}\label{sears2}
\sum_{i=0}^d q^{i(i+f-g)}\qbins{a-i}{b}\qbins{c-i}{d-i}\qbins{e}{i}
\qbins{f}{g-i}
=\sum_{i=0}^d q^{i(i+c-d-g)}\qbins{a-i}{b-i}\qbins{c-g}{d-i}
\qbins{a-e}{i}\qbins{d+f-i}{g-i},
\end{equation}
for $a,b,c,d,e,f,g\in\Z$ such that $b+c=d+e+f$.
The proof of this requires the substitutions
$n\to d$, $a\to q^{-g}$, $b\to q^{-b}$, $c\to q^{e-a}$, 
$d\to q^{-a}$, $e\to q^{-d-f}$ and $f\to q^{c-g-d+1}$ in
\eqref{Searstrafo} and some simplifications.
As a consequence of \eqref{sears2}
\begin{equation*}
\text{LHS}\eqref{ink}
=\sum_{i,k}q^{\frac{1}{2}b^2+k(k+b)+i(i+a-k)}
\qbins{L+M-b-k-i}{L-i}\qbins{a}{k-i}\qbins{M-a}{i}\qbins{L-a-i}{L-a-b-k-i}.
\end{equation*}
Shifting $k\to k+i$ and then renaming $i$ as $k$ and $k$ as $i$ yields
\begin{equation*}
\text{LHS}\eqref{ink}
=q^{\frac{1}{2}b^2}\sum_k q^{k(k+a+b)}\qbins{M-a}{k}
\sum_i q^{i(i+b+k)}\qbins{L+M-b-2k-i}{L-k}\qbins{a}{i}
\qbins{L-a-k}{L-a-b-2k-i}.
\end{equation*}
The sum over $i$ can be carried out by \eqref{qS} leading to the
right side of \eqref{ink}.
\end{proof}

\subsection{Proof of Theorem \ref{mainthm3}}\label{sec33}
As with the Theorems~\ref{mainthm} and \ref{mainthm2},
we first discuss the conditions imposed on the parameters.

On the left the summand vanishes unless $0\leq i\leq \min\{L-M,L\}$,
$i+2k\leq M$, $|b|\leq L-i-k$ and $|a|\leq k$.
Hence the left-hand side is nonzero if and only if $L\geq M$, $|a|\leq M/2$
and $|a|+|b|\leq L$.
By \eqref{rangeS} the right-hand side is nonzero if and only if
$|a+b|\leq L$, $|a|\leq M/2$ and $|a-b|\leq M$.
Because $|a|+|b|\leq L$ is equivalent to
$|a+b|\leq L$ and $|a-b|\leq L$,
both sides trivially vanish if any of the following three
conditions is violated: $|a|\leq M/2$, $|a+b|\leq L$ and 
$|a-b|\leq \max\{L,M\}$.
If these conditions are however satisfied, the mismatch between right and
left side needs to be repaired by imposing that $M\leq L$ and $|a-b|\leq M$.
(Note that $|a-b|\leq M$ implies $|a-b|\leq L$ when $M\leq L$.)

\begin{proof}[Proof of Theorem~\ref{mainthm3}]
After inserting the definitions \eqref{Sdef} and \eqref{Bdef} of $\S$ 
and $\B$, we use the symmetry \eqref{qbinsym} of the 
$q$-binomial coefficients to get
\begin{multline}\label{t3}
\sum_{i,k}q^{M i+k^2}\qbins{L-M}{L-M-i}\qbins{L+M-2i-2k}{M-i-2k}
\qbins{L-i-a+b}{k-a}\qbins{L-i+a-b}{L-i-k-b}\\
=q^{a^2}\sum_{k}
q^{k(k+a+b)}\qbins{L+M-a-b-2k}{L-a-b-2k}\qbins{M+a-b}{k}\qbins{M-a+b}{M-2a-k}.
\end{multline}
In view of the above discussion we may without loss of generality assume
that $M\leq L$ and $|a-b|\leq M$. This, together with the condition that 
the lower entries of all seven $q$-binomial coefficients are nonnegative,
implies that all the top-entries are nonnegative. Specifically, on the
left we have $L-M\geq 0$, $L+M-2i-2k=(M-i-2k)+(L-M-i)+M\geq 0$,
$L-i-a+b=(L-M-i)+(M-a+b)\geq 0$, $L-i+a-b=(L-M-i)+(M+a-b)\geq 0$.
Similarly, on the right we have $L+M-a-b-2k=(L-a-b-2k)+M\geq 0$, 
$M-a+b\geq 0$, $M-a+b\geq 0$. Consequently we can again assume definition 
\eqref{qbinmod} in the proof of \eqref{t3}.

By the simultaneous changes $i\to i+2k-M$ and $k\to M-i-k$
\begin{multline*}
\text{LHS}\eqref{t3} =
\sum_{i,k}q^{k^2+i(i+2k-M)} \\
\times\qbins{L-M}{L-i-2k}\qbins{L+M-2k}{i}
\qbins{L+M-i-2k-a+b}{M-i-k-a}\qbins{L+M-i-2k+a-b}{L-k-b}.
\end{multline*}
Transforming the sum over $i$ by \eqref{sears2} leads to
\begin{equation*}
\text{LHS}\eqref{t3}
=\sum_{i,k} q^{k^2+i(i+k+b)}\qbins{L+M-i-2k+a-b}{L-i-k-b}
\qbins{M-a+b}{M-i-k-a}\qbins{a-b}{i}\qbins{L-i-k-a}{L-i-2k},
\end{equation*}
which, by the simultaneous change $k\to i+a$ and $i\to k-i$, becomes
\begin{multline*}
\text{LHS}\eqref{t3}=q^{a^2}\sum_{k} q^{k(k+a+b)}\qbins{M-a+b}{M-k-2a} \\
\times \sum_i q^{i(i-k+a-b)}\qbins{L+M-i-k-a-b}{L-k-a-b}
\qbins{a-b}{k-i}\qbins{L-k-2a}{L-i-k-2a}.
\end{multline*}
An important difference between the definitions \eqref{qbin} and 
\eqref{qbinmod} of the $q$-binomial coefficients is that only the former 
satisfies the symmetry \eqref{qbinsym}. However, the modified $q$-binomials
do satisfy this symmetry provided $m+n\geq 0$. Now the above summand
vanishes if $M-k-2a<0$. Hence in the sum over $i$ we may assume
that $M-k-2a\geq 0$. This implies that $L-k-2a=(L-M)+(M-k-2a)\geq 0$,
so that we may replace $\qbins{L-k-2a}{L-i-k-2a}$ by $\qbins{L-k-2a}{i}$.
Then the sum over $i$ can be performed by \eqref{qS} 
resulting in the right-hand side of \eqref{t3}.
\end{proof}

\section{Two simple examples}\label{sec4}
\subsection{First example}
In our first application of the transformations \eqref{Tinv} and 
\eqref{TtoS} we start with the simplest possible identity for refined
$q$-trinomials.
\begin{lemma}\label{lemma21}
For $L,M\in\Z_{+}$
\begin{equation}\label{root1}
\sum_{j=-\infty}^{\infty}q^{j(j+1)}\{\T(L,M,2j,j)-\T(L,M,2j+2,j)\}=
\delta_{L,0}\delta_{M,0}.
\end{equation}
\end{lemma}
\begin{proof}
We begin with \cite[Lemma 3.1]{FLW97}
\begin{equation*}
\sum_{j=-\infty}^{\infty}\sum_{\tau=0}^1 (-1)^{\tau}
q^{j(j+1)}\B(L,M,j+\tau,j)=\delta_{L,0}\delta_{M,0}.
\end{equation*}
By \eqref{BtoT} this implies
\begin{equation*}
\sum_{j=-\infty}^{\infty}\sum_{\tau=0}^1 (-1)^{\tau}
q^{j(j+1)}\T(L,M,2j+2\tau,j)\!
=\sum_{\substack{n=0 \\n+L \text{ even}}}^M \!\!\!
q^{\frac{1}{2}n^2}\qbins{M}{n}\delta_{L,n}\delta_{M,0}
=\delta_{L,0}\delta_{M,0}.
\end{equation*}
\end{proof}
A single application of the transformation of Theorem~\ref{mainthm} yields
\begin{equation*}
\sum_{j=-\infty}^{\infty}
q^{\frac{1}{2}j(3j+2)}\{\T(L,M,3j,j)-\T(L,M,3j+2,j)\}=\delta_{L,0}.
\end{equation*}
Further iterating this, using \eqref{Tinv} and a simple induction argument,
shows that 
\begin{multline*}
\sum_{j=-\infty}^{\infty}q^{\frac{1}{2}j((k+3)j+2)}
\{\T(L,M,(k+3)j,j)-\T(L,M,(k+3)j+2,j)\} \\
=\sum_{r_1,\dots,r_{k-1}\geq 0}
q^{\frac{1}{2}(r_1^2+\cdots+r_k^2)}
\prod_{i=0}^{k-1}\qbins{L-r_{i+1}-\sum_{j=1}^{i-1}r_j}{r_i-r_{i+1}},
\end{multline*}
where $r_0:=M$, $r_k:=L-r_1-\cdots-r_{k-1}$ and $k\geq 1$.
Here and in the rest of the paper we adopt the convention that
$-\sum_{j=m}^{n-1}r_j=\sum_{j=n}^{m-1}r_j$ for $n<m$, so that
$$\qbins{L-r_1-\sum_{j=1}^{-1}r_j}{r_0-r_1}=\qbins{L+r_0-r_1}{r_0-r_1}
=\qbins{L+M-r_1}{M-r_1}.$$
We cannot turn the above polynomial identity into a nontrivial
$q$-series result because the large $L$ limit gives zero on either side.
What we can do is apply the transformation of Theorem~\ref{mainthm2}.
Note in particular that the condition $|a|\leq M$ if $|b|\leq M$ and 
$|a+b|\leq L$ does not pose a problem since 
$|(k+3)j|\geq |j|$ and $|(k+3)j+2|\geq |j|$ for $k\geq 0$.
By \eqref{TtoS} we thus find
\begin{multline}\label{LMBor}
\sum_{j=-\infty}^{\infty}\Bigl\{
q^{\frac{1}{2}j((k+3)(k+4)j+2)}\S(L,M,(k+4)j,(k+3)j) \\
-q^{\frac{1}{2}((k+3)j+2)((k+4)j+2)}\S(L,M,(k+4)j+2,(k+3)j+2)\Bigr\} \\
=\sum_{r_1,\dots,r_k\geq 0}
q^{\frac{1}{2}(r_1^2+\cdots+r_{k+1}^2)}
\qbins{L+M-r_1}{L}\qbins{L-r_2}{r_1}
\prod_{i=2}^k\qbins{r_1-r_{i+1}-\sum_{j=2}^{i-1}r_j}{r_i-r_{i+1}},
\end{multline}
where $r_{k+1}=r_1-r_2-\cdots-r_k$.
For $k=0$ the correct expression on the right side is $\qbins{L+M}{L}$.
It is also worth separately stating the $k=1$ case, namely
\begin{equation*}
\sum_{j=-\infty}^{\infty}(-1)^j 
q^{\frac{1}{2}j(5j+1)}\S(L,M,\lfloor (5j+1)/2 \rfloor,2j)
=\sum_{n=0}^M q^{n^2}\qbins{L+M-n}{L}\qbins{L-n}{n}.
\end{equation*}
This is the first of four doubly-bounded analogues of the first 
Rogers--Ramanujan identity that will be obtained in this paper.
Using \eqref{SlimL} to take the limit when $L$ tends to infinity
we find
\begin{equation}\label{Bress}
G(M,M;1,3/2,2)=\sum_{j=-\infty}^{\infty}(-1)^j q^{\frac{1}{2}j(5j+1)}
\qbins{2M}{M-2j}=\sum_{n=0}^M q^{n^2}\qbins{M}{n},
\end{equation}
a result originally due to Bressoud \cite[Eqs. (9)]{Bressoud81}.

If instead of $L$ we let $M$ become large and use \eqref{SlimM} we obtain
a result of Andrews \cite[Eq. (1.11)]{Andrews90}
\begin{equation*}
\sum_{j=-\infty}^{\infty}\bigl\{q^{j(10j+1)}\qbins{L}{5j}_2
-q^{(2j+1)(5j+2)}\qbins{L}{5j+2}_2\bigr\}
=\sum_{n=0}^{\infty} q^{n^2}\qbins{L-n}{n}.
\end{equation*}
For arbitrary $k$ equation \eqref{LMBor} is a generalization of
$q$-trinomial identities of \cite[Eq. (9.4)]{BMP98} and 
\cite[Prop. 4.5]{W99}.
Making the simultaneous replacements $r_i\to r_{i-1}-r_i$ ($i=2,\dots,k$) 
in \eqref{LMBor} and letting $L$ and $M$ tend to infinity gives 
identities for Virasoro characters of \cite[Cor. IV.1]{W99}.

Before we come to our next example let us point out that the above 
discussion can be repeated for the second Rogers--Ramanujan identity.
In particular one can show by a generalization of \eqref{TtoS} that
\begin{equation*}
\sum_{j=-\infty}^{\infty}(-1)^j 
q^{\frac{1}{2}j(5j+3)}\S(L,M,\lfloor 5j/2 \rfloor+1,2j+1)
=\sum_{n=0}^{M-1} q^{n(n+1)}\qbins{L+M-n-1}{L}\qbins{L-n-1}{n}.
\end{equation*}

\subsection{Second example}
Our next example uses a slightly more complicated-to-prove
identity as starting point.
\begin{lemma}
For $L,M\in\Z_{+}$
\begin{equation}\label{root2}
\sum_{j=-\infty}^{\infty}(-1)^j q^{\frac{1}{2}j(j+1)}\U(L,M,j,j)=\delta_{M,0}.
\end{equation}
\end{lemma}

\begin{proof}
Let $\sigma\in\{0,1\}$. Then for $L-\sigma/2\in\Z_+$ there holds
\begin{equation}\label{M0}
\sum_{j=-\infty}^{\infty}\sum_{\tau=0}^1 (-1)^{\tau} q^{j(2j+1)}
\B(L,M,j+\sigma/2,2j+\tau)=\delta_{M,0}.
\end{equation}
For $\sigma=0$ this is \cite[Cor. 3.2]{FLW97}. For $\sigma=1$ we
replace $L$ by $L+1/2$ so that we are to show that
\begin{equation*}
\sum_{j=-\infty}^{\infty}\sum_{\tau=0}^1 (-1)^{\tau} q^{j(2j+1)}
\qbins{L+M+j+\tau}{M+2j+\tau}\qbins{L+M-j-\tau+1}{M-2j-\tau}
=\delta_{M,0}.
\end{equation*}
By the $q$-binomial recurrence
\begin{equation}\label{qbinrec}
\qbins{m+n}{m}=\qbins{m+n-1}{m}+q^n\qbins{m+n-1}{m-1}
\end{equation}
the left side can be expanded as
\begin{equation*}
\sum_{j=-\infty}^{\infty}\sum_{\tau=0}^1 (-1)^{\tau} q^{j(2j+1)}
\qbins{L+M+j+\tau}{M+2j+\tau}
\Bigl\{\qbins{L+M-j-\tau}{M-2j-\tau}+q^{L+j+1}
\qbins{L+M-j-\tau}{M-2j-\tau-1}\Bigr\}.
\end{equation*}
The first term of the summand gives $\delta_{M,0}$ according to
\eqref{M0} with $\sigma=0$. The second term should thus vanish,
which readily follows by the variable changes $j\to -j-1$ and 
$\tau\to 1-\tau$.

If we now apply the transformation \eqref{BtoT} to \eqref{M0} we obtain
\begin{multline*}
\sum_{j=-\infty}^{\infty}\sum_{\tau=0}^1 (-1)^{\tau} q^{j(2j+1)}
\T(L,M,2j+\sigma,2j+\tau) \\
=\sum_{\substack{n=0 \\n+\sigma+L \text{ even}}}^M
q^{\frac{1}{2}n^2}\qbin{M}{n}\delta_{M,0}
=\delta_{M,0}\chi(L+\sigma\text{ even}),
\end{multline*}
where $\chi(\text{true})=1$ and $\chi(\text{false})=0$. 
Summing over $\sigma$, replacing $j\to (j-\tau)/2$ and using the 
definition \eqref{Udef} of $\U$ yields \eqref{root2}.
\end{proof}
By application of \eqref{Tinv} the identity \eqref{root2} transforms into 
\begin{equation*}
\sum_{j=-\infty}^{\infty}(-1)^j q^{\frac{1}{2}j(2j+1)}\U(L,M,2j,j)=
\qbin{L+M}{L}.
\end{equation*}
Letting $L$ tend to infinity using \eqref{UlimL} yields
a $q$-binomial identity equivalent to item H(2) in Slater's 
list of Bailey pairs \cite{Slater51}.
Next, by \eqref{Tinv} and induction
\begin{multline}\label{B1}
\sum_{j=-\infty}^{\infty}(-1)^j q^{\frac{1}{2}j((k+1)j+1)}\U(L,M,(k+1)j,j) \\
=\sum_{r_1,\dots,r_{k-1}\geq 0} q^{\frac{1}{2}(r_1^2+\cdots+r_{k-1}^2)}
\prod_{i=0}^{k-1}\qbins{L-r_{i+1}-\sum_{j=1}^{i-1}r_j}{r_i-r_{i+1}},
\end{multline}
for $k\geq 2$ and $r_0:=M$, $r_k:=0$. For $k=2$ and $L\to\infty$ this gives
a $q$-binomial identity of Rogers \cite{Rogers17} given as item B(1) 
in Slater's list. For $k=2,3$ and $M\to\infty$ this gives two
$q$-trinomial identities of Andrews \cite[Thm. 5.1]{Andrews90b},
\cite[Eq. (4.4)]{Andrews90}.
Before considering the identities arising when both $L$ and $M$ become
large, we will show that one can easily derive a variation of \eqref{B1}
using the following lemma.
\begin{lemma}
For $L,M\in\Z_{+}$
\begin{equation}\label{B2}
\sum_{j=-\infty}^{\infty}(-1)^j q^{3j(j+1)/2}\U(L,M,3j+1,j)=
q^M \sum_{r=0}^M q^{\frac{1}{2}r^2}
\qbins{L+M-r-1}{M-r}\qbins{L-1}{r}.
\end{equation}
\end{lemma}
For $L\to\infty$ this yields the Bailey pair B(2).
\begin{proof}
Fix $M$ and denote the left side of \eqref{B1} for $k=2$ by $f_L$
and the left side of \eqref{B2} by $g_L$. Using
\begin{equation}\label{qbinrec2}
\qbins{m+n}{m}=\qbins{m+n-1}{m-1}+q^m\qbins{m+n-1}{m}
\end{equation}
it then follows that 
\begin{align*}
f_{L-1}&=\sum_{j=-\infty}^{\infty}(-1)^j q^{\frac{1}{2}j(3j+1)}
\sum_{\tau=0}^1 \sum_{\substack{n=0 \\ n+L+j+\tau \text{ odd}}}^M 
q^{\frac{1}{2}n^2}\qbins{M}{n} \\
& \qquad \times \qbins{M+(L-j-n-\tau-1)/2}{M+j}
\qbins{M+(L+j-n+\tau-1)/2}{M-j} \\
&=q^{-M}\sum_{j=-\infty}^{\infty}(-1)^j q^{\frac{1}{2}j(3j+3)} 
\sum_{\tau=0}^1 \sum_{\substack{n=0 \\ n+L+j+\tau \text{ odd}}}^M 
q^{\frac{1}{2}n^2}\qbins{M}{n} \\ 
& \qquad \times \qbins{M+(L-j-n-\tau-1)/2}{M+j}
\Bigl\{\qbins{M+(L+j-n+\tau+1)/2}{M-j}-
\qbins{M+(L+j-n+\tau-1)/2}{M-j-1}\Bigr\}.
\end{align*}
The first term in the summand on the right yields $g_L$ and the second term
vanishes by the substitutions $j\to -j-1$ and $\tau\to 1-\tau$.
Hence $g_L=q^M f_{L-1}$. Since also the right sides of 
\eqref{B1} for $k=2$ and \eqref{B2} satisfy this equation we are done.
\end{proof}
We leave it to the reader to apply the transformation \eqref{Tinv}
to \eqref{B2} to obtain the variation of \eqref{B1} alluded to, but remark
that a single application of \eqref{Tinv} results in a
generalization of the $q$-trinomial identity 
\cite[Eq. (4.10)]{Andrews90}.

Next we take the large $L$ and $M$ limit in \eqref{B1}.
By \eqref{UlimLM} and the triple product identity \cite[Eq. (II.28)]{GR90}
\begin{equation}\label{tpi}
\sum_{k=-\infty}^{\infty} (-z)^k q^{\binom{k}{2}}=(z,q/z,q)_{\infty}
\end{equation}
this yields
\begin{equation*}
\sum_{m_1,\dots,m_{k-1}\geq 0} 
\frac{q^{\frac{1}{2}(M_1^2+\cdots+M_{k-1}^2)}}{(q)_{m_1}
\cdots (q)_{m_{k-1}}}
=\frac{(-q^{1/2})_{\infty}}{(q)_{\infty}}
(q^{k/2},q^{k/2+1},q^{k+1};q^{k+1})_{\infty},
\end{equation*}
with $M_i=m_i+\cdots+m_{k-1}$.
To connect this with more familiar $q$-series results we make use of
Lemma~\ref{lemred} to reduce the number of summation variables on the left.
First, for odd values of $k$ we use the expression for $f_{2k-1}(0)$
with $a=1$ given by the lemma to find the following theorem.
\begin{theorem}
For $k\geq 2$
\begin{align*}
\sum_{m_1,\dots,m_{2k-2}\geq 0}
\frac{q^{\frac{1}{2}(M_1^2+\cdots+M_{2k-2}^2)}}
{(q)_{m_1}\cdots (q)_{m_{2k-2}}}
&=\sum_{n_1,\dots,n_{k-1}\geq 0}
\frac{q^{N_1^2+\cdots+N_{k-1}^2}
(-q^{1/2-N_1})_{N_1}}{(q)_{n_1}\cdots(q)_{n_{k-1}}} \\
&=\prod_{\substack{j=1 \\ j\not\equiv 2\pmod{4} \\
j\not\equiv 0,\pm (2k-1) \pmod{4k}}}^{\infty}\frac{1}{(1-q^{j/2})}.
\end{align*}
\end{theorem}
The last two expressions of this theorem constitute Andrews' generalization
of the G\"ollnitz-Gordon identities \cite[Eq. (7.4.4)]{Andrews76}.
The equality of these with the first expression was conjectured by Melzer
\cite{Melzer94} and recently proven by Bressoud, Ismail and Stanton
\cite[Thm. 5.1; $i=k, a=1$]{BIS99} using different techniques.

Next, when $k$ is even we use the expression for $f_{2k}(0)$ with $a=1$.
\begin{theorem}
For $k\geq 1$
\begin{align*}
\sum_{m_1,\dots,m_{2k-1}\geq 0}
\frac{q^{\frac{1}{2}(M_1^2+\cdots+M_{2k-1}^2)}}
{(q)_{m_1}\cdots (q)_{m_{2k-1}}}
&=(-q^{1/2})_{\infty}\sum_{n_1,\dots,n_{k-1}\geq 0}
\frac{q^{N_1^2+\cdots+N_{k-1}^2}}{(q)_{n_1}\cdots(q)_{n_{k-1}}} \\
&=(-q^{1/2})_{\infty}
\prod_{\substack{j=1 \\ j\not\equiv 0,\pm k \pmod{2k+1}}}^{\infty}
\frac{1}{(1-q^j)}.
\end{align*}
\end{theorem}
The equality between the second and third expression is 
the well-known (first) Andrews--Gordon identity \cite{Andrews74}.

Returning to \eqref{B1} we apply the transformation \eqref{TtoS},
use definition \eqref{Vdef} and replace $k$ by $k-1$ to find
\begin{multline}\label{I3}
\sum_{j=-\infty}^{\infty}(-1)^j q^{\frac{1}{2}j(k(k+1)j+1)}
\V(L,M,(k+1)j,kj)\\
=\sum_{r_1,\dots,r_{k-1}\geq 0} q^{\frac{1}{2}(r_1^2+\cdots+r_{k-1}^2)}
\qbins{L+M-r_1}{L}\qbins{L-r_2}{r_1}\prod_{i=2}^{k-1}
\qbins{r_1-r_{i+1}-\sum_{j=2}^{i-1}r_j}{r_i-r_{i+1}},
\end{multline}
with $k\geq 2$ and $r_k:=0$.
Although the right sides of \eqref{B1} and \eqref{I3}
coincide for $k=2$, the large $L$ limit of \eqref{I3} for $k=2$
does not reproduce the Bailey pair B(1), but yields the 
pair I(3) due to Slater \cite{Slater51}.  
The large $M$ limit of this same identity yields \cite[Cor. 5.2]{Andrews90b}
of Andrews.
If for general $k$ we apply \eqref{VlimLM}
and collect even and odd powers of $q^{1/2}$ we obtain the 
Virasoro-character identity
\begin{equation*}
\sum_{\substack{r_1,\dots,r_{k-1}\geq 0 \\ 
\sigma+\sum_i r_i\text{ even}}}
\frac{q^{\frac{1}{2}(r_1^2+\cdots+r_{k-1}^2-\sigma)}}{(q)_{r_1}}
\prod_{i=2}^{k-1}\qbins{r_1-r_{i+1}-\sum_{j=2}^{i-1}r_j}{r_i-r_{i+1}}
=\begin{cases}
\chi^{(k,2k+2)}_{(k-1)/2,k+2\sigma}(q)
& \text{$k$ odd} \\[2mm]
\chi^{(k+1,2k)}_{k/2,k+2\sigma-1}(q)
& \text{$k$ even,}
\end{cases}
\end{equation*}
where $\sigma\in\{0,1\}$ and
\begin{equation*}
\chi_{r,s}^{(p,p')}(q)=\frac{1}{(q)_{\infty}}
\sum_{j=-\infty}^{\infty} 
\Bigl\{q^{j(pp'j+p'r-ps)}-q^{(pj+s)(p'j+r)}\Bigr\}.
\end{equation*}
This generalizes the identities (83) and (86) in Slater's list 
\cite{Slater52} of identities of the Rogers--Ramanujan type.

\section{Two not-so-simple examples}\label{sec5}
In the two examples of the previous section the initial
identities \eqref{root1} and \eqref{root2} were both 
straightforward consequences of known $q$-binomial identities. 
In this section we will give two further applications
of Theorems~\ref{mainthm} and \ref{mainthm2} that show that
not all irreducible identities for refined $q$-trinomials are trivial.
Apart from leading to more examples of identities of the
Rogers--Ramanujan type, this will result in the following 
remarkable pair of Virasoro-character identities: 
\begin{multline*}
\frac{1}{(q;q^2)_{\infty}}\Bigl\{
\chi^{(2k-1,2k+1)}_{k,k+1}(q^2)+q^{1/2}
\chi^{(2k-1,2k+1)}_{k,k}(q^2)\Bigl\} \\
=\sum_{r_1,\dots,r_{k-1}\geq 0}
\frac{q^{\frac{1}{2}(r_1^2+\cdots+r_{k-1}^2)}}
{(q)_{r_1}}\Bigl(\prod_{i=2}^{k-2}
\qbins{r_1-r_{i+1}-\sum_{j=2}^{i-1}r_j}{r_i-r_{i+1}}\Bigr)
\qbins{r_1+r_{k-1}-\sum_{j=2}^{k-2}r_j}{2r_{k-1}}_{q^{1/2}}
\end{multline*}
and
\begin{multline*}
\frac{1}{(-q^{1/2};q^{1/2})_{\infty}}\Bigl\{
\chi^{(2k-1,2k+1)}_{k,k+1}(q^{1/2})+q^{1/2}
\chi^{(2k-1,2k+1)}_{k,k-1}(q^{1/2})\Bigr\} \\
=\sum_{r_1,\dots,r_{k-1}\geq 0}
\frac{q^{\frac{1}{2}(r_1^2+\cdots+r_{k-1}^2)}}{(q)_{r_1}}
\Bigl(\prod_{i=2}^{k-2}
\qbins{r_1-r_{i+1}-\sum_{j=2}^{i-1}r_j}{r_i-r_{i+1}}\Bigl)
\qbins{\lfloor (r_1+r_{k-1}-\sum_{j=2}^{k-2}r_j)/2\rfloor}{r_{k-1}}_{q^2},
\end{multline*}
for $k\geq 3$ and $\lfloor x\rfloor$ the integer part of $x$.
The strange similarity between these two formulas and how the roles 
of $q^{1/2}$ and $q^2$ are interchanged in going from one to the other
is in our opinion quite amazing.

\subsection{A generalization of Bailey pair C(5)}
In our next example we take the following identity as starting point.
\begin{lemma}\label{lemc5}
For $L,M\in\Z_{+}$
\begin{equation}\label{root3}
\sum_{j=-\infty}^{\infty}(-1)^j q^{j(j+1)}\U(L,M,j,2j)=
q^{\binom{M}{2}}\qbin{L+M+1}{2M}_{q^{1/2}}.
\end{equation}
\end{lemma}
Taking the large $L$ limit using \eqref{UlimL} gives
\begin{equation*}
\sum_{j=-\infty}^{\infty}(-1)^jq^{j(j+1)}\qbin{2M}{M-2j}=
q^{\binom{M}{2}}(-q)_M.
\end{equation*}
This identity, which is equivalent to the Bailey pair C(5) 
of Rogers \cite{Rogers17,Slater51}, can be obtained as a specialization
of the nonterminating $q$-Dougall sum \cite[II.20]{GR90}.
This foreshadows that \eqref{root3} will not be as easy to prove
as our previous examples. 

\begin{proof}
By \eqref{qbin} and \eqref{rangeT} the lemma is trivially true for 
$M=0$, and in the following we may  assume $M\geq 1$. 
We now multiply both sides by $z^L$ and sum $L$ over 
the nonnegative integers. On the right this sum can be carried out thanks 
to \cite[(3.3.7)]{Andrews76}
\begin{equation}\label{qbt2}
\sum_{k=0}^{\infty} z^k \qbin{n+k}{k}=\frac{1}{(z)_{n+1}}.
\end{equation}
On the left we insert the definition \eqref{Udef} of $\U$ to obtain two terms.
In the second of these (corresponding to a sum over $\T(L,M,j+1,2j)$)
we change $j\to -j-1$ and use the symmetry \eqref{Qsymm}. 
Substituting the definition \eqref{Tdef}
of $\T$ and shifting $L\to 2L+n+j$ then gives
\begin{multline*}
\sum_{j=-\infty}^{\infty}\sum_{L=0}^{\infty}\sum_{\tau=0}^1\sum_{n=0}^M
(-1)^{j+\tau} z^{2L+j+n}q^{j(j+1)+\frac{1}{2}n^2}
\qbins{M}{n}\qbins{L+M+2j+2\tau}{L}\qbins{L+M-j-2\tau}{L+j} \\
=\frac{z^{M-1}q^{\binom{M}{2}}}{(z;q^{1/2})_{2M+1}}.
\end{multline*}
By \eqref{qbt} the sum over $n$ yields $(-zq^{1/2})_M$. Dividing by this 
term and using that $(-zq^{1/2})_M(z;q^{1/2})_{2M+1}=(z^2)_{2M+1}/(-z)_{M+1}$,
we are left to prove that
\begin{equation*}
\sum_{j=-\infty}^{\infty}\sum_{L=0}^{\infty}\sum_{\tau=0}^1
(-1)^{j+\tau} z^{2L+j}q^{j(j+1)}
\qbins{L+M+2j+2\tau}{L}\qbins{L+M-j-2\tau}{L+j}
=\frac{z^{M-1}q^{\binom{M}{2}}(-z)_{M+1}}{(z^2)_{2M+1}}.
\end{equation*}
We expand the right using \eqref{qbt} and \eqref{qbt2}
and equate coefficients of $z^a$. Renaming $L$ as $j$ and $M$ as $L$ 
this yields
\begin{multline*}
\sum_{j=0}^a \sum_{\tau=0}^1(-1)^{a+\tau}
q^{(2j-a)(2j-a-1)}
\qbins{L-3j+2a+2\tau}{j}\qbins{L+3j-a-2\tau}{a-j} \\
=q^{\binom{L}{2}} \sum_{i=0}^{\infty} q^{\binom{2i+L-a}{2}}
\qbins{2L+i}{i}\qbins{L+1}{2i+2L-a},
\end{multline*}
for $a\geq 0$ and $L\geq 1$.
Next we write $a=2M+\sigma$ for $\sigma\in\{0,1\}$,
and make the changes $j\to j+M+\sigma$ on the left and $i\to M-i$ 
on the right. Recalling \eqref{Brs} this gives
\begin{multline}\label{Bsigma}
\sum_{j=-\infty}^{\infty}\sum_{\tau=0}^1(-1)^{\tau}q^{2j(2j+1)}
\B_{0,\sigma}(L,M,4j+2\tau,j) \\
=q^{\binom{L}{2}} \sum_{i=0}^M q^{\binom{L-2i-\sigma}{2}}
\qbins{2L+M-i}{2L}\qbins{L+1}{2L-2i-\sigma},
\end{multline}
for $M\geq 0$ and $L\geq 1$.
To prove this identity we recall the following result due to
Gessel and Krattenthaler 
\cite[Thm. 12; $r=2$, $a=M$, $c=L$, $m=n=1=\epsilon=1$]{GK97}
\cite[Thm. 13; $r=2$, $a=M+1$, $c=L+1$, $m=\epsilon=1$, $n=0$]{GK97}
\begin{equation*}
\sum_{j=-\infty}^{\infty}\sum_{\tau=0}^1(-1)^{\tau}q^{j(3j+2-\sigma)}
\B_{\sigma,\sigma}(L,M,3j+2\tau,j)
=q^{L(L+\sigma-1)}\qbins{L+M+\sigma+1}{2L+\sigma},
\end{equation*}
for $L,M\geq 0$.
Applying Theorem~\ref{thmBurge3} with $r=s=\sigma$ yields \eqref{Bsigma}.
\end{proof}

If we now apply Theorem~\ref{mainthm} to Lemma~\ref{lemc5}
and simplify the resulting right-side by \eqref{A1} we obtain
\begin{equation}\label{C1}
\sum_{j=-\infty}^{\infty}(-1)^j q^{j(3j+1)}\U(L,M,3j,2j)=
\qbin{L+2M}{L}_{q^{1/2}},
\end{equation}
which in the large $L$ limit corresponds to the Bailey pair C(1)
\cite{Rogers17,Slater51}. We note that \eqref{C1} can be used to
also obtain a generalization of the Bailey pair C(2), or, 
more precisely, of a linear combination of C(1) and C(2).
\begin{lemma}\label{lemC2}
For $L,M\in\Z_{+}$
\begin{equation}\label{C2}
\sum_{j=-\infty}^{\infty}(-1)^j q^{j(3j+1)}\U(L,M,3j,2j+1)=
\qbin{L+2M-1}{L}_{q^{1/2}}.
\end{equation}
\end{lemma}
\begin{proof}
We will be rather brief, omitting some details.
Fixing $M$, let $f_L$ and $g_L$ denote the left side of \eqref{C1}
and \eqref{C2}, respectively. Then
\begin{multline*}
g_L=\sum_{j=-\infty}^{\infty}(-1)^j q^{j(3j+1)} \\ \times \sum_{\tau=0}^1
\sum_{\substack{n=0 \\ n+L+j+\tau \text{ even}}}^M 
q^{\frac{1}{2}n^2}\qbins{M}{n}
\qbins{M+(L+j-n-\tau+2)/2}{M+2j+1}\qbins{M+(L-j-n+\tau-2)/2}{M-2j-1}.
\end{multline*}
Applying \eqref{qbinrec2} to the second $q$-binomial yields two triple-sums.
One of these vanishes as it changes sign under the 
substitutions $j\to -j-1$ and $\tau\to 1-\tau$. Hence
\begin{multline*}
g_L=\sum_{j=-\infty}^{\infty}(-1)^j q^{j(3j+1)} \\ \times
\sum_{\tau=0}^1 \sum_{\substack{n=0 \\ n+L+j+\tau \text{ even}}}^M 
q^{\frac{1}{2}n^2}\qbins{M}{n}
\qbins{M+(L+j-n-\tau)/2}{M+2j}\qbins{M+(L-j-n+\tau-2)/2}{M-2j-1}.
\end{multline*}
Again by \eqref{qbinrec2} we can subtract this in a straightforward manner
from $f_L$. If in what then results we make the changes $j\to -j$ 
and $\tau\to 1-\tau$ we find $q^M f_{L-1}$. Hence $g_L=f_L-q^M f_{L-1}$.
If we now replace $f_L$ in this equation by the right side of \eqref{C1}
and use the recurrence \eqref{qbinrec} we find that $g_L$ equals
the right side of \eqref{C2}.
\end{proof}

Returning to \eqref{C1} it follows by \eqref{Tinv} and induction 
that for $k\geq 1$
\begin{multline}\label{C1p}
\sum_{j=-\infty}^{\infty}(-1)^j q^{j((2k+1)j+1)}\U(L,M,(2k+1)j,2j) \\
=\sum_{r_1,\dots,r_{k-1}\geq 0} q^{\frac{1}{2}\sum_{i=1}^{k-1}r_i^2}
\qbins{L+r_{k-1}-\sum_{j=1}^{k-2}r_j}{2r_{k-1}}_{q^{1/2}}
\prod_{i=0}^{k-2}
\qbins{L-r_{i+1}-\sum_{j=1}^{i-1}r_j}{r_i-r_{i+1}},
\end{multline}
with $r_0:=M$.
Letting $M$ tend to infinity we obtain a chain of
$q$-trinomial identities, the simplest of which
is equivalent to \cite[Thm. 6.1; (6.10)]{Andrews90b} by Andrews.
(The other identity in Andrews' theorem follows by a single iteration of
\eqref{C2} and again taking the large $M$ limit.)
When also $L$ tends to infinity we arrive at
\begin{equation*}
\sum_{n_1,\dots,n_{k-1}\geq 0} 
\frac{q^{N_1^2+\cdots+N_{k-1}^2}}
{(q^2;q^2)_{n_1}\cdots (q^2;q^2)_{n_{k-1}}(q;q^2)_{n_{k-1}}}
=\prod_{\substack{j=1 \\ j\not\equiv 2\pmod{4} \\
j\not\equiv 0,\pm 4k \pmod{8k+4}}}^{\infty}\frac{1}{(1-q^j)},
\end{equation*}
with $N_i=n_i+\cdots+n_{k-1}$.
Carrying out the same calculations starting with Lemma~\ref{lemC2}
leads to the analogous result
\begin{equation*}
\sum_{n_1,\dots,n_{k-1}\geq 0} 
\frac{q^{N_1^2+\cdots+N_{k-1}^2+2N_1+\cdots+2N_{k-1}}}
{(q^2;q^2)_{n_1}\cdots (q^2;q^2)_{n_{k-1}}(q;q^2)_{n_{k-1}+1}}
=\prod_{\substack{j=1 \\ j\not\equiv 2\pmod{4} \\
j\not\equiv 0,\pm 4 \pmod{8k+4}}}^{\infty}\frac{1}{(1-q^j)}.
\end{equation*}
For $k=2$ the above two identities are Rogers'
\cite[p. 330, Eq. (3)]{Rogers17} listed as items (79) and (96) in Slater's
list \cite{Slater52} of Rogers--Ramanujan identities.

When $k$ is odd we can apply Lemma~\ref{lemred} to obtain
the following equivalent pair of Rogers--Ramanujan identities
\begin{equation*}
\sum_{n_1,\dots,n_{k-1}\geq 0} 
\frac{q^{N_1^2+\cdots+N_{k-1}^2}}
{(q)_{n_1}\cdots (q)_{n_{k-1}}(q;q^2)_{n_{k-1}}}
=\prod_{\substack{j=1 \\ j\not\equiv 0,\pm 4k \pmod{8k-2}}}^{\infty}
\frac{1}{(1-q^j)}
\end{equation*}
and
\begin{equation*}
\sum_{n_1,\dots,n_{k-1}\geq 0} 
\frac{q^{N_1^2+\cdots+N_{k-1}^2+2N_1+\cdots+2N_{k-1}}}
{(q)_{n_1}\cdots (q)_{n_{k-1}}(q;q^2)_{n_{k-1}+1}}
=\prod_{\substack{j=1 \\ j\not\equiv 0,\pm 2 \pmod{8k-2}}}^{\infty}
\frac{1}{(1-q^j)}.
\end{equation*}
For $k=2$ these are \cite[p. 329, Eqs. (1.1) and (1.3)]{Rogers17} of Rogers,
corresponding to items (61) and (59) in Slater's list. We leave it to the
reader to carry out the corresponding rewritings when $k$ is even.

So far, the results of this section are nothing unusual; typical examples
of identities of the Rogers--Ramanujan type have been derived.
More exciting $q$-series results arise if we
apply \eqref{TtoS} to \eqref{C1p}. Replacing 
$k$ by $k-1$ this yields
\begin{multline}\label{C1pt}
\sum_{j=-\infty}^{\infty}(-1)^j q^{\frac{1}{2}j((4k^2-1)j+2)}
\V(L,M,(2k+1)j,(2k-1)j) \\
=\sum
q^{\frac{1}{2}\sum_{i=1}^{k-1} r_i^2}
\qbins{L+M-r_1}{L}\qbins{L-r_2}{r_1}
\qbins{r_1+r_{k-1}-\sum_{j=2}^{k-2}r_j}{2r_{k-1}}_{q^{1/2}}
\prod_{i=2}^{k-2}\qbins{r_1-r_{i+1}-\sum_{j=2}^{i-1}r_j}{r_i-r_{i+1}},
\end{multline}
where $k\geq 2$ and where the sum on the right is over 
$r_1,\dots,r_{k-1}\geq 0$.
Being of special interest, we separately state the $k=2$ case
\begin{equation}\label{RR1V}
\sum_{j=-\infty}^{\infty}(-1)^j q^{j(15j+2)}\V(L,M,5j,3j;q^2)
=\sum_{n=0}^M q^{n^2}\qbins{L+M-n}{L}_{q^2}\qbins{2L-n}{n}.
\end{equation}
Applying \eqref{TtoS} to \eqref{C2} we also have
\begin{equation}\label{RR1Vp}
\sum_{j=-\infty}^{\infty}(-1)^j q^{j(15j+2)}\V(L,M,5j+1,3j;q^2)
=\sum_{n=0}^M q^{n^2}\qbins{L+M-n}{L}_{q^2}\qbins{2L-n-1}{n},
\end{equation}
for $M\geq 1$.
These two formulas provide our second and third doubly-bounded analogue of 
the first Rogers--Ramanujan identity. 
If we follow Schur \cite{Schur17} and define the polynomials $e_n$
recursively as $e_n=e_{n-1}+q^{n-2} e_{n-2}$ with $e_0=0$ and $e_1=1$,
then the $M\to \infty$ limits of \eqref{RR1V} and \eqref{RR1Vp}
can be written as
\begin{equation}\label{eS}
e_{2L+\tau}=\sum_{j=-\infty}^{\infty}(-1)^j q^{j(15j+2)}
\Bigl\{\qbins{L}{5j+1-\tau}_{2;q^2}
+q^{6j+1}\qbins{L}{5j+2-\tau}_{2;q^2}\Big\}.
\end{equation}
where $\tau\in\{0,1\}$.
For $\tau=1$ this is equivalent to 
\cite[Eqs. (6.16) and (6.17)]{Andrews90b} of Andrews,
who remarked that on the right one can easily discern
even and odd powers of $q$. Even powers follow from $j$ even (odd) in the
first (second) term and odd powers follow from $j$ odd (even) in the
first (second) term. From the well-known combinatorial interpretation
of the Schur polynomial $e_n$ as the generating function of partitions
with difference between parts at least two and largest part not exceeding
$n-2$, it follows that the right-hand side of \eqref{eS} 
can be dissected to give the generating function of partitions of $m$ with
the parity of $m$ fixed, difference between parts at least two and 
largest part at most $2L+\tau-2$. It seems an interesting combinatorial 
problem to also interpret the identities \eqref{RR1V} and \eqref{RR1Vp}
in terms of restricted Rogers--Ramanujan partitions.

Letting not $M$ but $L$ tend to infinity in \eqref{RR1V} and \eqref{RR1Vp}
leads to the somewhat unusual Rogers--Ramanujan polynomial identity
\begin{equation*}
\sum_{j=-\infty}^{\infty}(-1)^j q^{j(15j+2)}
\Bigl\{\qbins{2M}{M-3j}_{q^2}+q^{6j+1}\qbins{2M}{M-3j-1}_{q^2}\Bigr\}
=\sum_{n=0}^M q^{n^2}\qbins{M}{n}(-q^{M-n+1})_n.
\end{equation*}

A slight modification of the previous derivations leads to analogous
results for the second Rogers--Ramanujan identity. To avoid repeating
ourselves we have chosen to only state the counterparts of \eqref{RR1V} 
and \eqref{RR1Vp}, given by
\begin{multline}\label{dgen}
\sum_{j=-\infty}^{\infty}(-1)^j q^{j(15j+4)}
\bigl\{\S(L,M,5j+\tau,3j;q^2)+q^{10j+3}\S(L,M,5j+3-\tau,3j+1;q^2)\bigr\}\\
=\sum_{n=0}^M q^{n(n+1)}\qbins{L+M-n}{L}_{q^2}\qbins{2L+\tau-n-2}{n},
\end{multline}
for $\tau\in\{0,1\}$, $L\geq 1-\tau$ and $M\geq 2-\tau$.
This can be proven using 
\begin{multline*}
\sum_{j=-\infty}^{\infty}(-1)^j q^{j(3j+2)}
\bigl\{\T(L,M,3j,2j+\tau)-\T(L,M,3j+2,2j+\tau)\bigr\}\\
=q^{L/2}\qbins{L+2M+\tau-2}{L}_{q^{1/2}},
\end{multline*}
for $L+M\geq \tau$, which follows from \eqref{C1} 
and \eqref{C2} using the recurrences \eqref{qbinrec} and \eqref{qbinrec2}.
When $M$ tends to infinity in \eqref{dgen} one obtains $q$-trinomial
identities for $d_{2L+\tau}$ where $d_n$ is again a
Schur polynomial, this time defined by $d_n=d_{n-1}+q^{n-2} 
d_{n-2}$ with $d_0=1$ and $d_1=0$. For $\tau=0$ these $q$-trinomial 
identities are equivalent to \cite[Eqs. (6.18) and (6.19)]{Andrews90b}.

After this intermezzo on polynomial analogues of the
Rogers--Ramanujan identities we return
to the general result \eqref{C1pt}, send $L$ and $M$ to infinity
and repeat the exercise of separating even and odd powers of $q^{1/2}$.
This yields the first identity stated in the introduction of this section.

\subsection{A generalization of Bailey pair G(4)}
Our final example before we come to the generalized Borwein conjecture
takes the following identity as starting point.
\begin{lemma}
For $L,M\in\Z_{+}$
\begin{equation}\label{root4}
\sum_{j=-\infty}^{\infty}(-1)^j q^{\frac{1}{4}j(j+1)}
\U(L,M,\lfloor\tfrac{1}{2}j\rfloor,j)
=(-1)^M q^{\frac{1}{2}M^2}\qbin{\lfloor (L+M+1)/2\rfloor}{M}_{q^2}.
\end{equation}
\end{lemma}
When $L$ becomes large this leads to
\begin{equation*}
\sum_{j=-\infty}^{\infty}(-1)^j q^{\frac{1}{4}j(j+1)}\qbin{2M}{M-j}
=(-1)^M q^{\frac{1}{2}M^2}(q^{1/2})_M,
\end{equation*}
which follows from the terminating $q$-Dougall sum \cite[(II.21)]{GR90}
and is equivalent to the Bailey pair G(4) of Rogers \cite{Rogers17,Slater51}.

\begin{proof}
We proceed as in the proof of Lemma~\ref{lemc5} and, assuming $M\geq 1$,
compute the generating function of \eqref{root4}.
Using the relation $(-zq^{1/2})_M(z^2;q^2)_{M+1}=
(z^2)_{2M+1}/(zq^{1/2})_M$ this yields to the identity
\begin{multline*}
\sum_{j=-\infty}^{\infty}\sum_{L=0}^{\infty}\sum_{\tau=0}^1(-1)^{\tau}
z^{2L+j}q^{\frac{1}{2}j(2j+1)} \Bigl\{
z\qbins{L+M+2j+\tau}{L}+\qbins{L+M+2j+\tau-1}{L-1}\Bigr\}
\qbins{L+M-j-\tau}{L+j} \\
=\frac{(-1)^M z^M(1+z)q^{\frac{1}{2}M^2}(zq^{1/2})_M}{(z^2)_{2M+1}}.
\end{multline*}
To the second term on the left we add the trivial identity
\begin{equation*}
\sum_{j=-\infty}^{\infty}\sum_{L=0}^{\infty}\sum_{\tau=0}^1(-1)^{\tau}
z^{2L+j}q^{L+\frac{1}{2}j(2j+1)} 
\qbins{L+M+2j+\tau-1}{L}\qbins{L+M-j-\tau}{L+j}=0.
\end{equation*}
(To prove this make the replacements $\tau\to 1-\tau$ and $L\to L-j$ 
followed by $j\to -j$.)
By the recurrence \eqref{qbinrec2} the second term then becomes equal 
to the first but without the $z$, so that both sides of the identity
can be divided by $1+z$. 
Equating coefficients of $z^a$ and renaming $L$ as $j$ and $M$ as $L$ 
then gives
\begin{multline*}
\sum_{j=0}^a \sum_{\tau=0}^1(-1)^{a+\tau}
q^{\frac{1}{2}(2j-a)(4j-2a-1)}
\qbins{L-3j+2a+\tau}{j}\qbins{L+3j-a-\tau}{a-j} \\
=q^{\frac{1}{2}L^2}\sum_{i=0}^{\infty}q^{\frac{1}{2}(2i+L-a)^2}
\qbins{2L+i}{i}\qbins{L}{2i+2L-a},
\end{multline*}
for $a\geq 0$ and $L\geq 1$. 
Setting $a=2M+\sigma$ for $\sigma\in\{0,1\}$ and changing $j\to j+M+\sigma$
on the left and $i\to M-i$ on the right results in
\begin{multline*}
\sum_{j=-\infty}^{\infty}\sum_{\tau=0}^1(-1)^{\tau}q^{j(4j+2\sigma+1)}
\B_{0,\sigma}(L,M,4j+\sigma+\tau,j) \\
=q^{\frac{1}{2}L^2} \sum_{i=0}^M q^{\frac{1}{2}(2i-L+\sigma)^2}
\qbins{2L+M-i}{2L}\qbins{L}{2L-2i-\sigma}.
\end{multline*}
This identity is simple consequence of Theorem~\ref{thmBurge3}
applied to
\begin{equation*}
\sum_{j=-\infty}^{\infty}\sum_{\tau=0}^1(-1)^{\tau}q^{j(3j+\sigma+1)}
\B_{\sigma,\sigma}(L,M,3j+\sigma+\tau,j)
=q^{L(L-\sigma)}\qbins{L+M+\sigma}{2L+\sigma},
\end{equation*}
which for $\sigma=0$ is due to Burge \cite[p. 217]{Burge93}
and for $\sigma=1$ to Gessel and Krattenthaler
\cite[Thm. 9; $r=2$, $a=M$, $c=L$, $m=\epsilon=1$, $n=0$]{GK97}.
\end{proof}

Applying the transformation \eqref{Tinv} and carrying out the resulting
sum on the right using \eqref{A2} leads to
\begin{equation}\label{G1}
\sum_{j=-\infty}^{\infty}(-1)^j q^{\frac{1}{4}j(3j+1)}
\U(L,M,\lfloor\tfrac{3}{2}j\rfloor,j)=\qbin{\lfloor L/2\rfloor+M}{M}_{q^2},
\end{equation}
which is a refinement of the Bailey pair G(1) \cite{Rogers17,Slater51}.
We remark that a calculation very similar to the one that led to
Lemma~\ref{lemC2} shows that \eqref{G1} implies a generalization of
the Bailey pair G(3);
\begin{equation*}
\sum_{j=-\infty}^{\infty}(-1)^j q^{\frac{3}{4}j(j+1)}
\U(L,M,3\lfloor\tfrac{1}{2}(j+1)\rfloor,j)
=q^M\qbin{\lfloor L/2\rfloor+M}{M}_{q^2}.
\end{equation*}
We will not pursue the consequences of this identity but content ourselves
with iterating just \eqref{G1}. By \eqref{Tinv} and induction
this gives
\begin{multline}\label{G}
\sum_{j=-\infty}^{\infty}(-1)^j q^{\frac{1}{4}j((2k+1)j+1)}
\U(L,M,\lfloor\tfrac{1}{2}(2k+1)j\rfloor,j) \\
=\sum_{r_1,\dots,r_{k-1}\geq 0}q^{\frac{1}{2}\sum_{i=1}^{k-1}r_i^2}
\qbins{\lfloor (L+r_{k-1}-\sum_{j=1}^{k-2}r_j)/2\rfloor}{r_{k-1}}_{q^2}
\prod_{i=0}^{k-2}\qbins{L-r_{i+1}-\sum_{j=1}^{i-1}r_j}{r_i-r_{i+1}},
\end{multline}
where $r_0:=M$ and $k\geq 2$.
We immediately consider the case when both $L$ and $M$ become large.
To shorten some of the resulting equations we also replace $k$
by $k+1$. Thanks to the triple product identity \eqref{tpi} this yields
\begin{equation}\label{S20}
\sum_{n_1,\dots,n_k\geq 0}
\frac{q^{\frac{1}{2}(N_1^2+\cdots+N_k^2)}}
{(q)_{n_1}\cdots (q)_{n_{k-1}}(q^2;q^2)_{n_k}}=
\frac{(q^{k/2+1/2},q^{k/2+1},q^{k+3/2};q^{k+3/2})_{\infty}}
{(-q)_{\infty}(q^{1/2};q^{1/2})_{\infty}},
\end{equation}
with $N_i=n_i+\cdots+n_k$. For $k=1$ this is \cite[page 330]{Rogers94},
given as entry (20) in Slater's list.

Once more we utilize Lemma~\ref{lemred}. When $k$ in \eqref{S20}
is even this implies 
\begin{equation*}
\sum_{n_1,\dots,n_k\geq 0} \frac{q^{N_1^2+\cdots+N_k^2}}
{(q)_{n_1}\cdots (q)_{n_k}(-q^{1/2};q^{1/2})_{2n_k}} 
=\frac{(q^{k+1/2},q^{k+1},q^{2k+3/2};q^{2k+3/2})_{\infty}}{(q)_{\infty}}.
\end{equation*}
For $k=1$ this is the first Rogers--Selberg identity
\cite[p. 339, Eq. (6.1)]{Rogers17}, \cite[Eq. (29)]{Selberg36} 
(or item (33) in Slater's list), and for 
general $k$ it is Paule's generalized Rogers--Selberg identity
\cite[Eq. (45); $r\to k+1$]{Paule85}.
Similarly, when $k$ is odd in \eqref{S20} we get
\begin{equation*}
\sum_{n_1,\dots,n_k\geq 0}
\frac{q^{N_1^2+\cdots+N_k^2}(-q^{1/2-N_1})_{N_1}}
{(q)_{n_1}\cdots (q)_{n_k} (-q^{1/2};q^{1/2})_{2n_k}}
=\frac{(q^k,q^{k+1/2},q^{2k+1/2};q^{2k+1/2})_{\infty}}
{(q^{1/2})_{\infty}(q^2;q^2)_{\infty}}.
\end{equation*}

Finally, applying \eqref{TtoS} to \eqref{G} and replacing $k$ by $k-1$ yields
\begin{multline*}
\sum_{j=-\infty}^{\infty}\bigl\{q^{\frac{1}{2}j((2k-1)(2k+1)j+1)}
\V(L,M,(2k+1)j,(2k-1)j) \\
-q^{\frac{1}{2}((2k-1)j+k-1)((2k+1)j+k)}
\V(L,M,(2k+1)j+k,(2k-1)j+k-1)\bigr\} \\
=\sum
q^{\frac{1}{2}\sum_{i=1}^{k-1}r_i^2}\qbins{L+M-r_1}{L}
\qbins{L-r_2}{r_1}
\qbins{\lfloor (r_1+r_{k-1}-\sum_{j=2}^{k-2}r_j)/2\rfloor}{r_{k-1}}_{q^2}
\prod_{i=2}^{k-2}\qbins{r_1-r_{i+1}-\sum_{j=2}^{i-1}r_j}{r_i-r_{i+1}},
\end{multline*}
where $k\geq 3$ and where the sum on the right is over 
$r_1,\dots,r_{k-1}\geq 0$. When $k=2$ the correct right side reads
$\sum_{r\geq 0}q^{\frac{1}{2}r^2}\qbins{L+M-r}{L}
\qbins{\lfloor L-r/2\rfloor}{L-r}_{q^2}$.
The second identity claimed in the introduction of this section follows
in the large $L$ and $M$ limit.

\section{The generalized Borwein conjecture}\label{sec6}
Finally we have come to our main application of the transformations of
section~\ref{sec2}, leading proofs of Theorems~\ref{T4}--\ref{T6}.

We begin by recalling some notation and results of \cite{W00b}.
Assume that $a,b$ are coprime integers such that $0<b<a$,
and define a nonnegative integer $n$ and positive integers $a_0,\dots,a_n$
as the order and partial quotients of the continued fraction representation
of $(a/b-1)^{\sign(a-2b)}$ ($\sign(0)=0$), i.e.,
\begin{equation*}
\Bigl(\frac{a}{b}-1\Bigr)^{\sign(a-2b)}=[a_0,\dots,a_n]=
a_0+\frac{1}{a_1+\cfrac{1}{\cdots+\cfrac{1}{a_n}}}.
\end{equation*}
For simplicity, $a_n$ will be fixed by requiring that $a_n\geq 2$ for 
$(a,b)\neq (2,1)$. (This is by no means necessary, see \cite{W00b}.)
We denote the continued fraction corresponding to $a,b$ by $\cf(a,b)$,
and note the obvious symmetry $\cf(a,b)=\cf(a,a-b)$.
We further define the partial sums $t_j=\sum_{k=0}^{j-1}a_k$ 
for $j=1,\dots,n+1$ and introduce $t_0=0$ and $d(a,b)=t_{n+1}$.
Finally we need a $d(a,b)\times d(a,b)$ matrix $\I(a,b)$ with entries
\begin{equation*}
\I(a,b)_{j,k}=\begin{cases}
\delta_{j,k+1}+\delta_{j,k-1} & \text{for $j\neq t_i$} \\
\delta_{j,k+1}+\delta_{j,k}-\delta_{j,k-1} & \text{for $j=t_i$}
\end{cases}
\end{equation*}
and a corresponding Cartan-type matrix $C(a,b)=2I-\I(a,b)$ where $I$ is
the $d(a,b)\times d(a,b)$ identity matrix.
Note that the matrix $\I(a,b)$ has the following block-structure:
\begin{equation*}
\I(a,b)=\left(\begin{array}{ccc|ccc|ccc}
&&&&&&&& \\
&T_{a_0}&&&&&&& \\
&&&-1&&&&& \\ \hline
&&1&&&&&& \\
&&&&\ddots&&&& \\
&&&&&&-1&& \\ \hline
&&&&&1&&& \\
&&&&&&&T_{a_n}& \\
&&&&&&&&
\end{array}\right)
\end{equation*}
where $T_i$ is the incidence matrix of the tadpole graph with $i$ vertices,
i.e., $(T_i)_{j,k}=\delta_{j,k+1} +\delta_{j,k-1}+\delta_{j,k}\delta_{j,i}$.

With the above notation we define a polynomial for each pair of
coprime integers $a,b$ such that $0<2b<a$ as follows:
\begin{equation}\label{Fab}
F_{a,b}(L,M)=\sum_{m\in\Z_+^{d(a,b)}}q^{L(L-2m_1)+m C(a,b)m}
\qbins{M+m_1-m_2}{L}\qbins{m_1-m_2}{n_1}
\prod_{j=2}^{d(a,b)}\qbins{\tau_j m_j+n_j}{\tau_j m_j}
\end{equation}
for $2b<a<3b$, and
\begin{equation}\label{Fab2}
F_{a,b}(L,M)=\sum_{m\in\Z_+^{d(a,b)}}q^{L(L-2m_1)+m C(a,b)m}
\qbins{M+m_2}{L}\qbins{m_2}{n_1}
\prod_{j=2}^{d(a,b)}\qbins{\tau_j m_j+n_j}{\tau_j m_j}
\end{equation}
for $0<3b\leq a$.
Here
\begin{equation*}
mC(a,b)m=\sum_{j,k=1}^{d(a,b)}m_j C(a,b)_{j,k}m_k=\sum_{j=0}^n
\Bigl(m_{t_j+1}^2+\sum_{k=t_j+1}^{t_{j+1}-1}(m_{k}-m_{k+1})^2\Bigr)
\end{equation*}
and $\tau_j=\tau_j(a,b)=2-\delta_{j,d(a,b)}$.
The auxiliary variables $n_j$ in the summand are integers
defined by the $(m,n)$-system
\begin{equation}\label{mn}
n_j=L\delta_{j,1}-\sum_{k=1}^{d(a,b)}C(a,b)_{j,k} m_k
\quad \text{for $j=1,\dots,d(a,b)$.}
\end{equation}

We are now prepared to state our first results of this section.
\begin{theorem}\label{T7}
For $L,M\in\Z_{+}$ and $a,b$ coprime integers such that $0<2b<a$,
\begin{equation}\label{eqab}
\sum_{j=-\infty}^{\infty}(-1)^j q^{\frac{1}{2}j((2ab+1)j+1)}\S(L,M,aj,2bj)
=F_{a,b}(L,M).
\end{equation}
\end{theorem}
A proof of this will be given in the next section.

In order to turn \eqref{eqab} into identities for $q$-binomial coefficients
we take the large $L$ limit. On the left this limit is easily computed
with the aid of \eqref{SlimL}. On the right some rewritings need
to be carried out first to cancel the term $L(L-2m_1)$ in the
exponent of $q$ in \eqref{Fab} and \eqref{Fab2}.
To this end we eliminate $m_1,\dots,m_{a_0}$ in favour of 
$n_1,\dots,n_{a_0}$.
By \eqref{mn} this yields
\begin{align}
m_j&=L-jm_{a_0+1}-\sum_{k=1}^{a_0} \min(j,k) n_k,
\qquad  j=1,\dots,a_0 \notag \\
n_{a_0+1}&=L-a_0 m_{a_0+1}-\sum_{k=1}^{a_0} k n_k
-\sum_{k=a_0+1}^{d(a,b)}C(a,b)_{a_0+1,k}m_k  \notag \\
n_j&=-\sum_{k=1}^{d(a,b)}C(a,b)_{j,k} m_k,
\qquad  j=a_0+2,\dots,d(a,b), \label{na0}
\end{align}
from which it follows that
\begin{equation*}
L(L-2m_1)+m C(a,b)m=\sum_{j=1}^{a_0} (N_j+m_{a_0+1})^2+
\sum_{j,k=a_0+1}^{d(a,b)} m_j C(a,b)_{j,k} m_k,
\end{equation*}
where $N_j=n_j+\cdots+n_{a_0}$.
Taking the large $L$ limit is now straightforward and if we define
$F_{a,b}(M)=(q)_M \lim_{L\to\infty}F_{a,b}(L,M)$ then
\begin{multline}\label{FabM}
F_{a,b}(M)
=\sum_{\substack{n_1,\dots,n_{a_0} \geq 0 
\\ m_{a_0+1},\dots,m_{d(a,b)}\geq 0}} \\ \times
\frac{q^{\sum_{j=1}^{a_0}(N_j+m_{a_0+1})^2
+\sum_{j,k=a_0+1}^{d(a,b)} m_j C(a,b)_{j,k} m_k}(q)_M}
{(q)_{M-2m_{a_0+1}-N_1-N_2}(q)_{n_1}\cdots(q)_{n_{a_0}}(q)_{2m_{a_0+1}}}
\prod_{j=a_0+2}^{d(a,b)}\qbins{\tau_j m_j+n_j}{\tau_j m_j}
\end{multline}
for $a>2b$. Here the auxiliary variables $n_j$ for $j\geq a_0+2$ are
given by \eqref{na0}. 
When $2b<a<3b$ there holds $a_0=1$ in which case $N_1=n_1$ and $N_2=0$. 
When $b=1$ there holds $d(a,1)=a_0=a-1$ in which case $m_{a_0+1}=0$.

Having defined the polynomial $F_{a,b}(M)$ we can now state the 
identities obtained when $L$ tends to infinity in Theorem~\ref{T7}.
\begin{corollary}\label{cor3}
For $L,M\in\Z_{+}$ and $a,b$ coprime integers such that $0<2b<a$,
\begin{equation*}
G(M,M;a/2,(a+1/b)/2,2b)=
\sum_{j=-\infty}^{\infty}(-1)^j q^{\frac{1}{2}j((2ab+1)j+1)}
\qbins{2M}{M-2bj}=F_{a,b}(M).
\end{equation*}
\end{corollary}
From \eqref{FabM} it follows that for $2b<a<3b$ the above right-hand side 
can be written as
\begin{equation*}
\sum_{n_1,m_2,\dots,m_{d(a,b)}\geq 0}
q^{(n_1+m_2)^2+\sum_{j,k=2}^{d(a,b)} m_j C(a,b)_{j,k} m_k}
\qbins{M}{n_1}\qbins{M-n_1}{2m_2}
\prod_{j=3}^{d(a,b)}\qbins{\tau_j m_j+n_j}{\tau_j m_j}\geq 0,
\end{equation*}
leading to our next corollary.
\begin{corollary}\label{cor4}
$G(M,M;b/2,(b+1/a)/2,2a)\geq 0$ for $a,b$ coprime integers such 
that $2a<b<3a$.
\end{corollary}
The reason for interchanging $a$ and $b$ in comparison with
Corollary~\ref{cor3} will become clear shortly.

To derive further positivity results form Theorem~\ref{T7}
we replace $q$ by $1/q$ in \eqref{eqab} using \eqref{Qdual}.
Defining the polynomial $f_{a,b}(L,M)$ by \eqref{Fab}
and \eqref{Fab2} but with $mC(a,b)m$ replaced by $\bar{m}C(a,b)m$
with $\bar{m}=(m_1,\dots,m_{d(a,b)-1},0)$
(so that $\bar{m}C(a,b)m=mC(a,b)m+m_{d(a,b)}(m_{d(a,b)-1}-m_{d(a,b)})$),
the following result arises.
\begin{corollary}
For $L,M\in\Z_{+}$ and $a,b$ coprime integers such that $0<2b<a$,
\begin{equation}\label{eqab2}
\sum_{j=-\infty}^{\infty}(-1)^j q^{\frac{1}{2}j((2ab-1)j+1)}
\S(L,M,aj,2bj)=f_{a,b}(L,M).
\end{equation}
\end{corollary}
The large $L$ limit can be taken following the same steps as
before, and if we define $f_{a,b}(M)=(q)_M \lim_{L\to\infty}f_{a,b}(L,M)$
then $f_{a,b}(M)$ for $b\neq 1$ is given by \eqref{FabM} with
$\sum_{j,k=a_0+1}^{d(a,b)} m_j C(a,b)_{j,k} m_k$
replaced by $\sum_{j,k=a_0+1}^{d(a,b)} \bar{m}_j C(a,b)_{j,k} m_k$
where $\bar{m}_j=m_j(1-\delta_{j,d(a,b)})$.
This leads to the following analogue of Corollary~\ref{cor3}.
\begin{corollary}\label{cor5}
For $L,M\in\Z_{+}$ and $a,b$ coprime integers such that $1<2b<a$,
\begin{equation*}
G(M,M;(a-1/b)/2,a/2,2b)=
\sum_{j=-\infty}^{\infty}(-1)^j q^{\frac{1}{2}j((2ab-1)j+1)}
\qbins{2M}{M-2bj}=f_{a,b}(M).
\end{equation*}
\end{corollary}
For $2b<a<3b$ the above right-hand side can be written as
\begin{equation*}
\sum_{n_1,m_2,\dots,m_{d(a,b)}\geq 0}
q^{(n_1+m_2)^2+\sum_{j,k=2}^{d(a,b)} \bar{m}_j C(a,b)_{j,k} m_k}
\qbins{M}{n_1}\qbins{M-n_1}{2m_2}
\prod_{j=3}^{d(a,b)}\qbins{\tau_j m_j+n_j}{\tau_j m_j}\geq 0.
\end{equation*}
Hence $G(M,M;(a-1/b)/2,a/2,2b)\geq 0$ for $2b<a<3b$.
If we apply the symmetry \eqref{Gsym} followed by the duality
\begin{equation*}
G(M,M;\alpha,\beta,K;1/q)=q^{-M^2}G(M,M;K-\alpha,K-\beta,K;q)
\end{equation*}
and make the simultaneous replacements $a\to 4a-b$ and $b\to a$,
we obtain the following result.
\begin{corollary}\label{cor6}
$G(M,M;b/2,(b+1/a)/2,2a)\geq 0$ for $a,b$ coprime integers such 
that $a<b<2a$.
\end{corollary}

In Corollary~\ref{cor5} $b=1$ is excluded, requiring a different
treatment. Namely, if we express the summand of $f_{k+1,1}(L,M)$ 
in terms of $n_1,\dots,n_k$ we find that \eqref{eqab2} becomes 
\begin{multline*}
\sum_{j=-\infty}^{\infty}(-1)^j q^{\frac{1}{2}j((2k+1)j+1)}
\S(L,M,(k+1)j,2j) \\=\sum_{n_1,\dots,n_k\geq 0}
q^{N_1^2+\cdots+N_k^2+n_k\tilde{N}_k}
\qbins{M+\tilde{N}_2}{L}\qbins{\tilde{N}_2}{n_1} 
\qbins{n_k+\tilde{N}_k}{n_k}
\prod_{j=2}^{k-1}\qbins{n_j+2\tilde{N}_j}{n_j},
\end{multline*}
with $N_j=n_j+\cdots+n_k$, $\tilde{N}_j=L-N_1-\cdots-N_j$ and $k\geq 2$.
Noting that in the large $L$ limit only the
terms with $n_k=0$ contribute to the sum on the right,
this can be recognized as a doubly-bounded analogue of the
(first) Andrews--Gordon identity, with $k=2$ corresponding
to our fourth doubly-bounded version of the first
Rogers--Ramanujan identity. 
Taking $L\to\infty$ yields
\begin{equation*}
\sum_{j=-\infty}^{\infty}(-1)^j q^{\frac{1}{2}j((2k+1)j+1)}
\qbins{2M}{M-2j}=\sum_{n_1,\dots,n_{k-1}\geq 0}
\frac{q^{N_1^2+\cdots+N_{k-1}^2}(q)_M}
{(q)_{M-N_1-N_2}(q)_{n_1}\cdots(q)_{n_{k-1}}},
\end{equation*}
which for $k=2$ is Bressoud's identity \eqref{Bress}.
We remark that it follows from \eqref{FabM} that also Theorem~\ref{T7}
for $b=1$ is a doubly-bounded analogue of the (first) Andrews--Gordon
identity, with the same large $L$ limit as above, but with $k=2$ excluded.

Corollaries~\ref{cor4} and \ref{cor6} are the $\bar{a}=(2)$ and
$\bar{a}=(0,2)$ instances of Theorem~\ref{T6}.
We can use the full set of transformations of section~\ref{sec2}
to derive polynomial identities that imply all of Theorem~\ref{T6}.
However, the notation required to describe these more general
identities is a lot more complicated than the already involved 
notation needed in the definitions of $F_{L,M}$, $F_M$, $f_{L,M}$ and $f_M$.
Since we trust that the previous examples (and the proof of
Theorem~\ref{T7} given in the next section) illustrate how one
can, in principle, obtain explicit representations for all of the
polynomials occurring in Theorem~\ref{T6} we will in the following
take a shortcut and prove the remainder of the theorem without
first deriving explicit polynomial identities.

We begin by noting that the Corollaries~\ref{cor1} and \ref{cor2}
imply the following simple lemma.
\begin{lemma}\label{lemGtoG}
If $G(M,M;\alpha,\beta,K)\geq 0$ then $G(M,M;\alpha',\beta',K')\geq 0$
with 
\begin{align}\label{GtoG1}
\alpha'&=\alpha/2+K,&\quad \beta'&=\beta/2+K, &\quad  K'&=2K \\
\intertext{or}
\alpha'&=(\alpha+K)/2,&\quad \beta'&=(\beta+K)/2, &\quad  K'&=2K.
\label{GtoG2}
\end{align}
\end{lemma}
\begin{proof}
By Corollary~\ref{cor1} and the assumption that 
$G(M,M;\alpha,\beta,K)\geq 0$,
\begin{align*}
0&\leq \sum_{k=0}^{\infty}C_{M,k}(q)G(k,k;\alpha,\beta,K) \\
&=\sum_{j=-\infty}^{\infty} (-1)^j 
q^{Kj((\alpha+\beta)j+\alpha-\beta)/2}
\sum_{k=0}^{\infty}C_{M,k}(q)\qbins{2k}{k-Kj} \\
&=\sum_{j=-\infty}^{\infty} (-1)^j 
q^{Kj((4K+\alpha+\beta)j+\alpha-\beta)/2}
\qbins{2M}{M-2Kj} 
=G(M,M;\alpha',\beta',K'),
\end{align*}
with $\alpha',\beta'$ and $K'$ given by \eqref{GtoG1}.
Using Corollary~\ref{cor2}
instead of \ref{cor1} and copying the above steps
leads to a proof of the second statement.
\end{proof}
Next we iterate Lemma~\ref{lemGtoG} to arrive at a binary tree of positivity
results.
\begin{proposition}\label{prop}
For $n\geq 0$, let $\bar{a}=(a_0,\dots,a_n)\in\Z^{n+1}$ such that
$a_0\geq 0$, $a_1,\dots,a_{n-1}\geq 1$ and $a_n\geq 2-\delta_{n,0}$.
Then $G(M,M;\alpha,\beta,K)\geq 0$ implies that 
$G(M,M;\alpha',\beta',K')\geq 0$ with
\begin{equation*}
\alpha'=(\alpha+\mu K)/2^{d-1},\quad
\beta'=(\beta+\mu K)/2^{d-1},\quad K'=2^{d-1}K
\end{equation*}
where $d,\mu\in\Z_+$ are given by \eqref{dmu}.
\end{proposition}
The reason for referring to this as a binary tree is that in the proof
given below $G(M,M;\alpha,\beta,K)$ corresponds to the ``initial condition''
$\bar{a}=(1)$ (for which $d=1$ and $\mu=0$), \eqref{GtoG1} corresponds to the
transformation $(a_0,a_1,\dots,a_n)\to (a_0+1,a_1,\dots,a_n)$ and 
\eqref{GtoG2} corresponds to $(a_0,a_1,\dots,a_n)\to 
(0,a_0+1,a_1,\dots,a_n)$.

Having said this, it is clear that we could equally well have
chosen a different labelling of the triples $(\alpha',\beta',K')$ in 
the tree. For example, we could have chosen to let \eqref{GtoG2} correspond 
to the transformation $(a_0,a_1,\dots,a_n)\to (a_0+1,a_1,\dots,a_n)$ and 
\eqref{GtoG1} to $(a_0,a_1,\dots,a_n) \to 
(0,a_0+1,a_1,\dots,a_n)$.
On the sequences $\bar{a}$ this corresponds to the involution 
$(1)\to (1)$ and 
\begin{equation*}
(a_0,a_1,\dots,a_n)\to
(0,1^{a_0-1},2,1^{a_1-2},2,\dots,2,1^{a_{n-1}-2},2,1^{a_n-3},2),
\end{equation*}
which leaves $d$ invariant.
Here $1^{\mu}$ stands for $\mu$ repeated ones,
and $$(0,1^{-1})^{\mu}(2,1^{-1})^{\nu},2$$
(which can occur for $\mu\in\{0,1\}$ and $\nu\geq 0$;
and stands for $0,1^{-1}$ 
repeated $\mu$ times followed by $2,1^{-1}$ repeated $\nu$ times
followed by a $2$) has to be identified with the single integer $\nu-\mu+2$.
For example, $(1,1,2)\to(0,1^0,(2,1^{-1})^2,2)=(0,4)$
and $(0,4)\to (0,1^{-1},2,1,2)=(1,1,2)$.
With this different choice of labelling the statement of the
Proposition remains the same except that the expression for $\mu$
would then be
\begin{equation*}
\mu=\frac{1}{3}(4^{d-1}-1)+\sum_{i=1}^n 4^{a_i+\cdots+a_n-2}.
\end{equation*}
We do however believe that no relabelling is possible that would
simplify the claim of the proposition.

\begin{proof}[Proof of Proposition~\ref{prop}]
Given a sequence $\bar{a}=(a_0,\dots,a_n)$ we write 
$d=a_0+\cdots+a_n=|\bar{a}|$. We will also make the dependence of
$\alpha',\beta'$ and $K'$ on $\bar{a}$ explicit by writing 
$\alpha'_{\bar{a}},\beta'_{\bar{a}}$ and $K'_{\bar{a}}$.
Similarly we will write $\alpha_{\bar{a}}$, $\beta_{\bar{a}}$ and 
$K_{\bar{a}}$ whenever necessary, but 
$\alpha_{(1)}$, $\beta_{(1)}$ and $K_{(1)}$ will always be just
$\alpha$, $\beta$ and $K$.

As already remarked above, for $\bar{a}=(1)$, which is the only admissible
sequence for which $|\bar{a}|=1$, the proposition is trivially true.
We will now proceed by induction on $d$ and assume that the proposition is
true for sequences $\bar{a}$ with $|\bar{a}|=d$.
Now there are two types of (admissible) sequences $\bar{a}'$ that have 
$|\bar{a}'|=d+1$. Either it is of the form $\bar{a}'=(a_0+1,a_1,\dots,a_n)$ 
or it is of the form $\bar{a}'=(0,a_0+1,a_1,\dots,a_n)$ where in both
cases $a_0\geq 0$. (If $n=0$ in the latter case, then $a_0\geq 1$.)

First assume $\bar{a}'=(a_0+1,a_1,\dots,a_n)$. Since $|\bar{a}'|=d+1$ 
the sequence $\bar{a}=(a_0,\dots,a_n)$ has $|\bar{a}|=d$ and by our 
induction hypothesis the proposition holds for this $\bar{a}$.
If we now apply the \eqref{GtoG1} case of Lemma~\ref{lemGtoG} to $\bar{a}$,
 and then use the induction hypothesis, we find the proposition to be 
true for
\begin{equation*}
K'=2K_{\bar{a}}=2(2^{d-1}K)=2^d K=K_{\bar{a}'}
\end{equation*}
and
\begin{equation*}
\alpha'=\alpha_{\bar{a}}/2+K_{\bar{a}}=
(\alpha+\mu_{\bar{a}}K)/2^d+2^{d-1}K=
(\alpha+\mu_{\bar{a}'}K)/2^d=\alpha_{\bar{a}'},
\end{equation*}
where we have used that $\mu_{\bar{a}'}-\mu_{\bar{a}}=4^{d-1}$.
Repeating the last calculation with $\alpha$ replaced by $\beta$
also shows that $\beta'=\beta_{\bar{a}'}$.

Next assume $\bar{a}'=(0,a_0+1,a_1,\dots,a_n)$. Since $|\bar{a}'|=d+1$
the sequence $\bar{a}=(a_0,\dots,a_n)$ has $|\bar{a}|=d$ and by our 
induction hypothesis the proposition holds for this $\bar{a}$.
If we now apply the \eqref{GtoG2} case of Lemma~\ref{lemGtoG}
to $\bar{a}$, and then use the induction hypothesis,
we find the proposition to be true for $K'=K_{\bar{a}'}$
and
\begin{equation*}
\alpha'=(\alpha_{\bar{a}}+K_{\bar{a}})/2=
(\alpha+\mu_{\bar{a}}K)/2^d+2^{d-2}K=
(\alpha+\mu_{\bar{a}'}K)/2^d=\alpha_{\bar{a}'},
\end{equation*}
where we have used that $\mu_{\bar{a}'}-\mu_{\bar{a}}=4^{d-1}$.
Repeating this with $\alpha$ replaced by $\beta$
shows that $\beta'=\beta_{\bar{a}'}$.
\end{proof}

Using the Theorems~\ref{T1}--\ref{T3} as input to Proposition~\ref{prop}
the Theorems~\ref{T4}--\ref{T6} easily follow, and we only
present the proof of Theorem~\ref{T6}.
\begin{proof}[Proof of Theorem~\ref{T6}]
By Theorem~\ref{T1} there holds 
$G(M,M,\tilde{b},\tilde{b}+1/\tilde{a},\tilde{a})\geq 0$ for 
$\tilde{a},\tilde{b}$ coprime integers such that $0<\tilde{b}<\tilde{a}$.
By Proposition~\ref{prop} it is therefore true that
$G(M,M,\alpha,\beta,K)\geq 0$ with
$\alpha=(\tilde{b}+\mu\tilde{a})/2^{d-1}$,
$\beta=(\tilde{b}+1/\tilde{a}+\mu\tilde{a})/2^{d-1}$ and
$K=2^{d-1}\tilde{a}$. Defining $a=\tilde{a}$ and $b=\tilde{b}+\mu\tilde{a}$
gives the statement of the theorem. (Since
$\tilde{a}$ and $\tilde{b}$ are coprime, so are $a$ and $b$.)
\end{proof}

A close scrutiny of the set of admissible sequences $\bar{a}$
allows for a slight reformulation of Theorem~\ref{T6}.
First note there are $2^{d-1}$ sequences with fixed $d$.
For example, when $d=4$ we have the following eight
sequences in reverse lexicographic order: $S_4:=\{(4),(2,2),(1,3),(1,1,2),
(0,4),(0,2,2),(0,1,3),(0,1,1,2)\}$. Now observe that these eight sequences
form four pairs, with a typical pair given by $(a_0,\dots,a_{n-1},a_n)$ and
$(a_0,\dots,a_{n-1},a_n-2,2)$ with $a_n\geq 3$. (The only exception is the
pair $(2)$ and $(0,2)$ for $d=2$, which corresponds to $a_n=a_0=2$.)
Moreover, if $b$ and $b'$ form such a pair, with $b>b'$ in reverse
lexicographical order, then $\mu_b=\mu_{b'}+1$ with $\mu_b\equiv 2\pmod{4}$.
For example, the eight values of $\mu$ corresponding to the elements
of the set $S_4$ are given by $42,41,38,37,26,25,22,21$.
We can thus reformulate the above theorem as follows.
\begin{theorem}\label{T8}
For $n\geq 0$, let $\bar{a}=(a_0,\dots,a_n)\in\Z^{n+1}$ such that
(ii) $a_0\geq 0$, $a_1,\dots,a_{n-1}\geq 1$ and
$a_n\geq 3-\delta_{n,0}$. Then 
\begin{equation*}
G(M,M;b/2^{d-1},(b+1/a)/2^{d-1},2^{d-1}a)\geq 0,
\end{equation*}
for $a,b$ coprime integers such that
$(\mu-1)a<b<(\mu+1)a$, with $d$ and $\mu$ given by \eqref{dmu}.
\end{theorem}
Note that in comparison with Theorem~\ref{T6} the case $\bar{a}=(1)$ is 
now excluded.
We also note that the discussion leading to the Theorem~\ref{T8}
does not quite justify the claim of the theorem.
After all, what we really have argued is that $a,b$ must satisfy
$(\mu-1)a<b<(\mu+1)a$ and $b\neq \mu a$.
The following proof is to show that this latter condition can be dropped.
\begin{proof}
Since $a,b$ are positive, coprime integers, the only
solution to $b=\mu a$ is given by $(a,b)=(1,\mu)$.
So the problem is to show that 
$G(M,M;\mu/2^{d-1},(\mu+1)/2^{d-1},2^{d-1})\geq 0$ for each admissible
sequence $\bar{a}$. Now fix $\bar{a}$.
From $\mu_{\bar{a}}\equiv 2\pmod{4}$
it follows that $\mu_{\bar{a}}/2$ is an odd integer.
Hence we can apply Theorem~\ref{T6} with $(a,b)=(2,\mu_{\bar{a}}/2)$
and sequence $\bar{a}'=(a_0,\dots,a_n-1)$.
Since $d_{\bar{a}'}=d_{\bar{a}}-1$ and $\mu_{\bar{a}'}=(\mu_{\bar{a}}-2)/4$
the inequality $\mu_{\bar{a}'} a<b<(\mu_{\bar{a}'}+1)a$ translates into 
$\mu_{\bar{a}}-2<\mu_{\bar{a}}<\mu_{\bar{a}}+2$ and is 
therefore satisfied, as required by Theorem~\ref{T6}. 
But with the above choice for $a,b$ and $\bar{a}'$
Theorem~\ref{T6} tells us that 
$G(M,M;\mu_{\bar{a}}/2^{d_{\bar{a}}-1},(\mu_{\bar{a}}+1)/2^{d_{\bar{a}}-1},
2^{d_{\bar{a}}-1})\geq 0$.
\end{proof}

\section{Proof of Theorem \ref{T7}}\label{sec7}
In addition to the definitions of the previous section
we set $d(1,1)=1/\cf(1,1)=0$.
To also facilitate computations involving continued fractions,
we sometimes, by abuse of notation, write $\cf(a,b)=(a/b-1)^{\sign(a-2b)}$.

For coprime integers $a,b$ such that $1\leq b\leq a$ define 
\begin{equation}\label{Gab}
G_{a,b}(L,M)=\sum_{m\in\Z_+^{d(a,b)}}q^{m C(a,b)m}
\qbins{\tau_0 L+M-m_1}{\tau_0 L}
\prod_{j=1}^{d(a,b)}\qbins{\tau_j m_j+n_j}{\tau_j m_j}
\end{equation}
for $1\leq b\leq a\leq 2b$ (so that $G_{1,1}(L,M)=\qbins{L+M}{L}$) and
\begin{equation}\label{Gaab}
G_{a,b}(L,M)=\sum_{m\in\Z_+^{d(a,b)}}q^{L(L-2m_1)+m C(a,b)m}
\qbins{L+M+m_1}{2L}
\prod_{j=1}^{d(a,b)}\qbins{\tau_j m_j+n_j}{\tau_j m_j}
\end{equation}
for $a\geq 2b$. 

Our proof of Theorems \ref{T7} relies on the following identity
for the polynomial $G_{a,b}$ \cite[Lem. 3.1 and Thm 3.1]{W00b}.
\begin{theorem}\label{T9}
For $L,M\in\Z_{+}$ and $a,b$ coprime integers such that $1\leq b\leq a$,
\begin{equation*}
\sum_{j=-\infty}^{\infty}(-1)^j q^{\frac{1}{2}j((2ab+1)j+1)}
\B(L,M,aj,bj)=G_{a,b}(L,M).
\end{equation*}
\end{theorem}

\subsection{Proof of Theorem~\ref{T7} for $2b<a<3b$}
Take Theorem~\ref{T9}, replace $a,b$ by $\bar{a},\bar{b}$
and apply \eqref{BtoT} followed by \eqref{TtoS}.
With the notation $a=2\bar{a}+\bar{b}$ and $b=\bar{a}$ this leads to
\begin{multline}\label{Sab}
\sum_{j=-\infty}^{\infty}(-1)^j q^{\frac{1}{2}j((2ab+1)j+1)}
\S(L,M,aj,2bj) \\
=\sum_{\substack{r_1,r_2\geq 0 \\ r_1+r_2 \text{ even}}}
q^{\frac{1}{2}(r_1^2+r_2^2)}\qbins{L+M-r_1}{L}\qbins{L-r_1}{r_2}
f_{\bar{a},\bar{b}}(\tfrac{1}{2}(r_1-r_2),L-r_1).
\end{multline}
Next insert the expression for $G_{\bar{a},\bar{b}}$ given in 
\eqref{Gab} and \eqref{Gaab}. First, when $1\leq \bar{b}<\bar{a}<2\bar{b}$
($\tau_0=2$ since $(\bar{a},\bar{b})\neq (1,1)$),
\begin{multline}\label{exp}
\text{RHS}\eqref{Sab}=
\sum_{\substack{r_1,r_2\geq 0 \\ r_1+r_2 \text{ even}}}
\sum_{m\in\Z_+^{d(\bar{a},\bar{b})}}
q^{\frac{1}{2}(r_1^2+r_2^2)+m C(\bar{a},\bar{b})m} \\
\times
\qbins{L+M-r_1}{L}\qbins{L-r_1}{r_2}
\qbins{L-r_2-m_1}{r_1-r_2}
\prod_{j=1}^{d(\bar{a},\bar{b})}
\qbins{\bar{\tau}_j m_j+n_j}{\bar{\tau}_j m_j},
\end{multline}
with $\bar{\tau}_j=\tau_j(\bar{a},\bar{b})$ and
$n_j$ given by \eqref{mn} with $L\to (r_1-r_2)/2$ and $a,b\to\bar{a},\bar{b}$.
Now relabel $m_j\to m_{j+2}$ and $n_j\to n_{j+2}$,
then replace $r_1\to L-m_1+m_2$ and $r_2\to L-m_1-m_2$, and introduce
the auxiliary variables $n_1=L-m_1-m_2$ and $n_2=m_1-m_2-m_3$.
Since $\bar{a}<2\bar{b}$ and $a>2b$ one finds
$\cf(a,b)=1+1/(1+1/\cf(\bar{a},\bar{b}))$ and thus
$\cf(a,b)=[1,1,\alpha_0,\dots,\alpha_n]$ with 
$\cf(\bar{a},\bar{b})=[\alpha_0,\dots,\alpha_n]$
(with $\alpha_n\geq 2$ since $(\bar{a},\bar{b})\neq (2,1)$).
This implies $d(a,b)=d(\bar{a},\bar{b})+2$,
$\tau_j=\tau_j(a,b)=\bar{\tau}_{j+2}$ and
\begin{equation*}
C(a,b)=\left(\begin{array}{rr|rrr}
1&1&&& \\ -1&1&1&& \\ \hline &-1&&& \\ &&&
C(\bar{a},\bar{b})& \\ &&&& \end{array}\right)
\end{equation*}
and therefore $\text{RHS}\eqref{Sab}=\text{RHS}\eqref{Fab}$.
Eliminating $\bar{a}$ and $\bar{b}$ in $1\leq \bar{b}<\bar{a}<2\bar{b}$
in favour of $a$ and $b$ yields the condition $5b/2<a<3b$.

Next, when $2\leq 2\bar{b}\leq \bar{a}$, one again finds \eqref{exp},
but with an additional $(r_1-r_2)(r_1-r_2-4m_1)/4$ in the exponent of $q$.
By the same variable changes as before this yields an extra $m_2(m_2-2m_3)$
in the exponent of $q$. Since $\bar{a}\geq 2\bar{b}$ and $a>2b$ one finds
$\cf(a,b)=1+1/(1+\cf(\bar{a},\bar{b}))$ and thus
$\cf(a,b)=[1,1+\alpha_0,\alpha_1,\dots,\alpha_n]$ with 
$\cf(\bar{a},\bar{b})=[\alpha_0,\dots,\alpha_n]$.
(For $(\bar{a},\bar{b})=(2,1)$ one finds $\cf(a,b)=\cf(5,2)=[1,2]$
which has $\alpha_n\geq 2$ as it should.)
Hence $d(a,b)=d(\bar{a},\bar{b})+2$, $\tau_j=\tau_j(a,b)=\bar{\tau}_{j+2}$
and
\begin{equation*}
C(a,b)=\left(\begin{array}{rr|rrr}
1&1&&& \\ -1&2&-1&& \\ \hline &-1&&& \\ &&&
C(\bar{a},\bar{b})& \\ &&&& \end{array}\right)
\end{equation*}
and thus again $\text{RHS}\eqref{Sab}=\text{RHS}\eqref{Fab}$.
This time, however, $2\leq 2\bar{b}\leq\bar{a}$ leads to the condition
$2b<a\leq 5b/2$.

\subsection{Proof of Theorem~\ref{T7} for $a\geq 3b$}
Take Theorem~\ref{T9} with $a,b$ replaced by $\bar{a},\bar{b}$,
use the symmetry \eqref{Bsymm} and then apply \eqref{BtoT} followed by 
\eqref{TtoS}.
With the notation $a=\bar{a}+2\bar{b}$ and $b=\bar{b}$ this gives
\begin{multline}\label{Sab2}
\sum_{j=-\infty}^{\infty}(-1)^j q^{\frac{1}{2}j((2ab+1)j+1)}
\S(L,M,aj,2bj) \\
=\sum_{\substack{r_1,r_2\geq 0 \\ r_1+r_2 \text{ even}}}
q^{\frac{1}{2}(r_1^2+r_2^2)}
\qbins{L+M-r_1}{L}\qbins{L-r_1}{r_2}
G_{\bar{a},\bar{b}}(L-r_1,\tfrac{1}{2}(r_1-r_2)).
\end{multline}

Now insert the explicit expression for $G_{\bar{a},\bar{b}}$.
First, when $1\leq \bar{b}\leq \bar{a}<2\bar{b}$,
\begin{multline}\label{exp2}
\text{RHS}\eqref{Sab2}=
\sum_{\substack{r_1,r_2\geq 0 \\ r_1+r_2 \text{ even}}}
\sum_{m\in\Z_+^{d(\bar{a},\bar{b})}}
q^{\frac{1}{2}(r_1^2+r_2^2)+m C(\bar{a},\bar{b})m} \\
\times \qbins{L+M-r_1}{L}\qbins{L-r_1}{r_2}
\qbins{\tau_0(L-r_1)+\frac{1}{2}(r_1-r_2)-m_1}{\tau_0(L-r_1)}
\prod_{j=1}^{d(\bar{a},\bar{b})}
\qbins{\bar{\tau}_j m_j+n_j}{\bar{\tau}_j m_j},
\end{multline}
with $n_j$ given by \eqref{mn} with $L\to L-r_1$ and 
$(a,b)\to (\bar{a},\bar{b})$.

The next step is to relabel $m_j\to m_{j+2}$ and $n_j\to n_{j+2}$,
then to replace $r_1\to L-m_2$ and $r_2\to L-2m_1+m_2$, and 
to introduce the auxiliary variables $n_1=L-2m_1+m_2$ and 
$n_2=m_1-m_2-m_3$.
Since $\bar{a}<2\bar{b}$ and $a>2b$ one finds
$\cf(a,b)=2+1/\cf(\bar{a},\bar{b}))$ and thus
$\cf(a,b)=[2,\alpha_0,\dots,\alpha_n]$ with
$\cf(\bar{a},\bar{b})=[\alpha_0,\dots,\alpha_n]$.
This implies $d(a,b)=d(\bar{a},\bar{b})+2$,
$\tau_j=\tau_j(a,b)=\bar{\tau}_{j+2}$ and
\begin{equation*}
C(a,b)=\left(\begin{array}{rr|rrr}
2&-1&&& \\ -1&1&1&& \\ \hline &-1&&& \\ &&&
C(\bar{a},\bar{b})& \\ &&&& \end{array}\right)
\end{equation*}
($C(3,1)=((2,-1),(-1,1))$)
and thus $\text{RHS}\eqref{Sab2}=\text{RHS}\eqref{Fab2}$.
Writing $1\leq \bar{b}\leq \bar{a}<2\bar{b}$ in terms of
$a$ and $b$ yields the condition $3b\leq a<4b$.

Next, when $\bar{a}\geq 2\bar{b}$, one again finds
\eqref{exp2}, but with an additional $(L-r_1)(L-r_1-2m_1)$ in the exponent
of $q$. By the same variable changes as above this yields
an extra $m_2(m_2-2m_3)$ in the exponent of $q$.
Since $\bar{a}\geq 2\bar{b}$ and $a>2b$ this yields
$\cf(a,b)=2+\cf(\bar{a},\bar{b}))$ and thus
$\cf(a,b)=[2+\alpha_0,\alpha_1,\dots,\alpha_n]$ with
$\cf(\bar{a},\bar{b})=[\alpha_0,\dots,\alpha_n]$.
This implies $d(a,b)=d(\bar{a},\bar{b})+2$,
$\tau_j=\tau_j(a,b)=\bar{\tau}_{j+2}$ and
\begin{equation*}
C(a,b)=\left(\begin{array}{rr|rrr}
2&-1&&& \\ -1&2&-1&& \\ \hline &-1&&& \\ &&&C(\bar{a},\bar{b})& \\ 
&&&& \end{array}\right),
\end{equation*}
and once again $\text{RHS}\eqref{Sab2}=\text{RHS}\eqref{Fab2}$.
The condition $\bar{a}\geq 2\bar{b}$ implies $a\geq 4b$.

\appendix 
\section{Some simple summation formulas}\label{secqhyper}
In this appendix some simple identities used in the main text are proven.

Our first result is nothing but a corollary of the $q$-Saalsch\"utz sum
\cite[Eq. (II.12)]{GR90}
\begin{equation}\label{qSaal}
\sum_{k=0}^n \frac{(a,b,q^{-n})_k\, q^k}{(q,c,abq^{1-n}/c)_k}=
\frac{(c/a,c/b)_n}{(c,c/ab)_n}.
\end{equation}
Specializing $n\to M$,
$a\to q^{-L/2}$, $b\to q^{-(L+1)/2}$ and $c\to q^{1/2}$, 
and making some simplifications, yields
\begin{equation}\label{A1}
\sum_{i=0}^M q^{\frac{1}{2}i(2i-1)}\qbin{L+M-i}{L}\qbin{L+1}{2i}_{q^{1/2}}=
\qbin{L+2M}{L}_{q^{1/2}}.
\end{equation}
Also the next identity follows from \eqref{qSaal}, albeit with a bit 
more effort,
\begin{equation}\label{A2}
\sum_{i=0}^M(-1)^i q^{i^2}
\qbin{L+M-i}{L}\qbin{\lfloor (L+1)/2\rfloor}{i}_{q^2}
=\qbin{\lfloor L/2\rfloor+M}{M}_{q^2}.
\end{equation}
When $L$ is even this follows from the substitutions
$n\to M$, $a=-b\to q^{-L/2}$ and $c\to -q$ in \eqref{qSaal}.
To obtain \eqref{A2} for $L$ odd we denote the left side of 
\eqref{A2} by $f_{L,M}$ and note that by \eqref{qbinrec2}
$f_{2L-1,M}=f_{2L,M}-q^{2L} f_{2L,M-1}$.
Since we already proved \eqref{A2} for even $L$ we may on the right
replace $f_{2L,M}$ and $f_{2L,M-1}$ using \eqref{A2}.
By \eqref{qbinrec} this yields $f_{2L-1,M}=\qbins{L+M-1}{M}_{q^2}$ 
completing the proof. We note that \eqref{A2} for odd $L$ can also be
viewed as a corollary of a basic hypergeometric summation, given
by the $\tau=2$ instance of
\begin{equation*}
{_3\phi_2}\Bigl[\genfrac{}{}{0pt}{}{a^{1/2},-a^{1/2},q^{-n}}
{-q,aq^{1-n}};q,q^{\tau}\Bigr]:=
\sum_{k=0}^n \frac{(a;q^2)_k (q^{-n})_k q^{\tau k}}
{(q^2;q^2)_k(aq^{1-n})_k}
=\frac{q^{(2-\tau)n}(1/a;q^2)_n}{(-q,1/a)_n}
\end{equation*}
true for $\tau\in\{0,1\}$. The proof of this almost balanced summation 
proceeds along the same lines as the proof of \eqref{A2}.

Simple as it is, our final summation formula ---
used in the main text to simplify multiple sums --- appears to be new.
We remark that it can also be used very effectively to reduce the number of
(independent) entries in Slater's list of 130 Rogers--Ramanujan-type
identities~\cite{Slater52}. 
For $M_i=m_i+\cdots+m_k$ define
\begin{equation*}
f_k(m_k)=\sum_{m_1,\dots,m_{k-1}\geq 0} 
\frac{a^{M_1+\cdots+M_k}q^{\frac{1}{2}(M_1^2+\cdots+M_k^2)}}
{(q)_{m_1}\cdots(q)_{m_{k-1}}}.
\end{equation*}
\begin{lemma}\label{lemred}
Let $N_i=n_i+\cdots+n_k$. Then
\begin{align*}
f_{2k-1}(n_k)&=\sum_{n_1,\dots,n_{k-1}\geq 0}
\frac{a^{2N_1+\cdots+2N_k}q^{N_1^2+\cdots+N_k^2}(-q^{1/2-N_1}/a)_{N_1}}
{(q)_{n_1}\cdots(q)_{n_{k-1}}(-aq^{1/2})_{n_k}} \\
f_{2k}(n_k)&=(-aq^{1/2})_{\infty}\sum_{n_1,\dots,n_{k-1}\geq 0} 
\frac{a^{2N_1+\cdots+2N_k}q^{N_1^2+\cdots+N_k^2}}
{(q)_{n_1}\cdots (q)_{n_{k-1}}(-aq^{1/2})_{n_k}}. 
\end{align*}
\end{lemma}
\begin{proof}
When $k$ is odd we replace $k$ by $2k-1$ and introduce
new summation variables $n_1,\dots,n_{k-1}$ and $t_1,\dots,t_{k-1}$ as
$n_i=m_{2i-1}+m_{2i}$ and $t_i=m_{2i}$.
Using the notation of the lemma this gives
\begin{align*}
f_{2k-1}(n_k)&=\sum_{n_1,\dots,n_{k-1}\geq 0}
a^{N_1}q^{\frac{1}{2}N_1^2}
\prod_{i=1}^{k-1} a^{2N_{i+1}}q^{N_{i+1}^2}\sum_{t_i=0}^{n_i} \frac{a^{t_i}
q^{\frac{1}{2}t_i(t_i+2N_{i+1})}} {(q)_{t_i}(q)_{n_i-t_i}} \\
&=\sum_{n_1,\dots,n_{k-1}\geq 0}
\frac{a^{N_1+2N_2+\cdots+2N_k}q^{\frac{1}{2}N_1^2+N_2^2+\cdots+N_k^2}}
{(q)_{n_1}\cdots(q)_{n_{k-1}}}\prod_{i=1}^{k-1} (-aq^{N_{i+1}+1/2})_{n_i} \\
&=\sum_{n_1,\dots,n_{k-1}\geq 0}
\frac{a^{2N_1+\cdots+2N_k}q^{N_1^2+\cdots+N_k^2}
(-q^{1/2-N_1}/a)_{N_1}}{(q)_{n_1}\cdots(q)_{n_{k-1}}(-aq^{1/2})_{n_k}},
\end{align*}
where the second equality follows from the $q$-binomial theorem \eqref{qbt}
with $x=-aq^{N_{i+1}+1/2}$ and the last equality follows from
\begin{align*}
\prod_{i=1}^{k-1}(-aq^{N_{i+1}+1/2})_{n_i}
&=\prod_{i=1}^{k-1}\frac{(-aq^{1/2})_{N_i}}{(-aq^{1/2})_{N_{i+1}}} \\
&=\frac{(-aq^{1/2})_{N_1}}{(-aq^{1/2})_{N_k}}
=a^{N_1}q^{\frac{1}{2}N_1^2}\frac{(-q^{1/2-N_1}/a)_{N_1}}
{(-aq^{1/2})_{n_k}}.
\end{align*}
Next, when $k$ is even, we replace $k$ by $2k$ and introduce
new variables $n_1,\dots,n_{k-1}$, $t_1,\dots,t_{k-1}$ and $s$
as $n_i=m_{2i}+m_{2i+1}$, $t_i=m_{2i-1}$ and $s=\sum_{j=1}^k m_{2j-1}$.
With the notation $T_i=s-t_1-\cdots-t_i$ and 
$t_k=T_{k-1}=s-t_1-\cdots-t_{k-1}$ this yields
\begin{equation*}
f_{2k}(n_k)=
\sum_{\substack{n_1,\dots,n_{k-1}\geq 0\\s,t_1,\dots,t_{k-1}\geq 0}}
\frac{a^{s+2N_1+\cdots+2N_k}q^{\frac{1}{2}s^2+sn_k+N_1^2+\cdots+N_k^2+
\sum_{i=1}^{k-1}t_i(N_i-n_k-T_i)}}
{(q)_{T_{k-1}}\prod_{i=1}^{k-1}(q)_{t_i}(q)_{n_i-t_{i+1}}}.
\end{equation*}
Now define 
\begin{equation*}
g(n_1,\dots,n_l,s)=\sum_{t_1,\dots,t_l\geq 0}\frac{1}{(q)_{T_l}}
\prod_{i=1}^l \frac{q^{t_i(n_i+\cdots+n_l-T_i)}}
{(q)_{t_i}(q)_{n_i-t_{i+1}}},
\end{equation*}
where $t_{l+1}=T_l$. Obviously, $g(s)=1/(q)_s$. By the 
$q$-Chu--Vandermonde sum \eqref{qCV}
with $a\to\infty$, $n\to T_{l-1}$ and $c\to q^{n_l-T_{l-1}+1}$
it follows that 
$g(n_1,\dots,n_l,s)=g(n_1,\dots,n_{l-1})/(q)_{n_l}$. Hence
\begin{equation*}
g(n_1,\dots,n_l,s)=\frac{1}{(q)_s(q)_{n_1}\cdots(q)_{n_l}}.
\end{equation*}
When $l=k-1$ we insert this in the expression for $f_{2k}(n_k)$ to get
\begin{equation*}
f_{2k}(n_k)=
\sum_{n_1,\dots,n_{k-1},s\geq 0}
\frac{a^{s+2N_1+\cdots+2N_k}q^{\frac{1}{2}s^2+sn_k+N_1^2+\cdots+N_k^2}}
{(q)_s(q)_{n_1}\cdots(q)_{n_l}}.
\end{equation*}
Performing the sum over $s$ by the $L\to\infty$ limit of
\eqref{qbt} settles the second claim of the lemma.
\end{proof}

\bibliographystyle{amsplain}

\end{document}